\documentclass[11pt, a4paper]{article}

\usepackage{amsmath,amsfonts,amssymb,amsthm}
\usepackage{graphicx}
\usepackage[utf8]{inputenc}
\usepackage[T1]{fontenc}
\usepackage{algorithm}
\usepackage{algorithmic}
\usepackage{floatpag}
\usepackage{float}
\usepackage{graphics}
\usepackage{xcolor}
\usepackage{listings}
\usepackage[final]{showlabels}
\usepackage{multirow}

\usepackage[nocomma]{optidef}
\usepackage[colorlinks=true,linkcolor=blue]{hyperref}

\usepackage{tikz,pgfplots} 
\pgfplotsset{compat=1.3}
\usetikzlibrary{external}
\graphicspath{ {./figs/} }

\usepackage{caption}
\usepackage[T1]{fontenc}
\usepackage[absolute,overlay]{textpos}
\usepackage{hyperref}
\usepackage{url}
\usepackage{xcolor}
\usepackage{authblk}
\usepackage{blindtext}
\usepackage{array,longtable}
\usepackage{subfig}
\allowdisplaybreaks 
\makeatletter
\g@addto@macro{\UrlBreaks}{\UrlOrds}
\makeatother
\usepackage{comment}
\usepackage{booktabs}

\let\realverbatim=\verbatim
\let\realendverbatim=\endverbatim
\renewcommand\verbatim{\par\addvspace{6pt plus 2pt minus 1pt}\realverbatim}
\renewcommand\endverbatim{\realendverbatim\addvspace{6pt plus 2pt minus 1pt}}

\newsavebox{\astrutbox}
\sbox{\astrutbox}{\rule[-5pt]{0pt}{20pt}}

\newcommand{\mbf}[1]{\mathbf{#1}}

\theoremstyle{plain}
\newtheorem{theorem}{Theorem}[section]   

\newtheorem*{corollary*}{Corollary}  

 \theoremstyle{definition}
 \newtheorem{remark}[theorem]{Remark}

 \newtheorem*{problem*}{Problem}    
 \newtheorem*{example*}{Example}

\usepackage{hyperref} \hypersetup{ colorlinks = true, citecolor = {blue},
linkcolor = {green} }
\usepackage[backend=biber, style=numeric, sortlocale=auto, 
doi=false, isbn=false, url=false,eprint=false,
giveninits=true,maxnames=6]{biblatex}
\title{Impact of environmental constraints\\
in hydrothermal energy planning}

\author{%
{Lu\'{i}s Felipe Bueno\thanks{Universidade Federal de São Paulo, Campus São José dos Campos, lfelipebueno@gmail.com}}, 
{Andr\'e Luiz Diniz\thanks{Centro de Pesquisas de Energia Elétrica}},\\ 
{Rafael Durbano Lobato\thanks{Câmara de Comercialização de Energia Elétrica}}, 
{Claudia Sagastizábal\thanks{Universidade Estadual de Campinas}}, 
{Kenny Vinente\thanks{Universidade Federal do Amazonas}}
}

\begin{document}

\label{firstpage}
\maketitle

\begin{footnotesize}

\textbf{Studygroup}: 
9th Brazilian Study Group with Industry, 6-10 March 2023, São Carlos
\url{http://www.cemeai.icmc.usp.br/WSMPI/}

\textbf{Communicated by:}
Francisco Louzada Neto and Jos\'e Alberto Cuminato

\textbf{Industrial partners:}
Center of Research on Electric Energy (CEPEL) and
Chamber of Electric Energy Commercialization (CCEE)

\textbf{Presenter:}
Academic: Luis Felipe Bueno, Kenny Vinente. 
Industrial: Andr\'e Diniz, Rafael Lobato

\textbf{Team:}
{
Gabriel Vinicius Bacci, Zeray Hagos Gebrezabher, Kevin Felipe Oliveira (Universidade de São Paulo~-- USP, Brazil), 
Robério da Rocha Barboza, Dr. Andr\'e Diniz (Centro de Pesquisas de Energia Elétrica -- CEPEL, Brazil),
Dr. Felipe Beltran, Gabriel Teixeira, Lucas Dagort (Norus, Brazil), 
		Dr. Lu\'{i}s Felipe Bueno, Humberto Gimenes, Lucas Julião, Maria Clara Couto Lorena, Dimary Moreno (Universidade Federal de São Paulo -- UNIFESP, Brazil), 
  Ruan Felipe Sousa (Universidade Federal do Rio de Janeiro -- UFRJ, Brazil), 
  Dr. Rafael Lobato (Câmara de Comercialização de Energia Elétrica -- CCEE, Brazil),	
  Renata Pedrini (Universidade Federal de Santa Catarina  -- UFSC, Brazil),
  Dr. Claudia Sagastizábal, Dr. Williams Yanez (Universidade Estadual de Campinas -- UNICAMP, Brazil),
  Dr. Kenny Vinente (Universidade Federal do Amazonas -- UFAM, Brazil). 
}

\textbf{Application}:
Energy/Utilities

\textbf{Tools:}
Mathematical optimization, Julia, JuMP, Gurobi



\end{footnotesize}

\abstract{
As a follow-up of the industrial problems dealt with in 2018, 2019, 2021 and 2022, 
in partnership with CCEE and  CEPEL, 
in 2023 the study group 
``Energy planning and environmental constraints'' focused on the
impact that prioritizing multiple uses of water has on the 
the electric energy production systems, specially in predominantly hydro systems, which is the case of Brazil. 

In order to model environmental constraints in the
long-term hydrothermal generation planning problem, the resulting large-scale multi-stage
linear programming problem was modelled in JuMP and solved by stochastic dual
dynamic programming. To assess if the development represented well the behavior
of the Brazilian power system, the Julia formulation first was benchmarked with Brazil's official model, Newave.
Environmental constraints were introduced in this problem by two different
approaches, one that represents the multiple uses of water by means of 0-1
variables, and another one that makes piecewise linear approximations of the
relevant constraints. Numerical results show that penalties of
slack variables strongly affect the obtained water values.
}

\section*{Nomenclature}
\textbf{Indices}
\begin{longtable*}{@{}p{1.0cm}p{11.5cm}@{}}
$t$ & index of periods (months). \\
$g$ & index of thermal plants.\\
$h$ & index of hydro plants.\\
$a$ & index of scenario.\\
$b(a)$ & index of the scenario that preceded $a$ in the previous time period.\\
$r(h)$ & index of the reference reservoir of plant $h$.\\
$\underline{z}$ & index of volume range used for minimum outflow constraints.\\ 
$\overline{z}$ & index of volume range used for maximum outflow constraints.\\ 
\end{longtable*}

\textbf{Variables}
\begin{longtable*}{@{}p{1.0cm}p{11.5cm}@{}}
$ph_{hta}$ & power output of hydro $h$ in period $t$  under scenario $a$ (MW). \\
$pt_{gta}$ & power output of thermal plant $g$ in period $t$  under scenario $a$ (MW). \\
$q_{hta}$ & turbine discharge of hydro $h$ in period $t$  under scenario $a$ (m$^3$/s). \\
$s_{hta}$ & spillage of hydro $h$ in period $t$  under scenario $a$ (m$^3$/s). \\
$v_{hta}$ & reservoir volume of hydro $h$ at period $t$ under scenario $a$ (hm$^3$). \\
$def_{ta}$ & power deficit  at period $t$  under scenario $a$ (MW). \\
$def_{hta}^{\min}$ & deficit  of the minimum outflow constraint for hydro $h$ at period $t$   under scenario $a$ (m$^3$/s). \\
$def_{hta}^{\max}$ & deficit  of the maximum outflow constraint for hydro $h$ at period $t$ under scenario $a$ (m$^3$/s). \\
$def_{hta}^{\rm Res}$ & deficit  of the maximum outflow constraint for hydro $h$ at period $t$ under scenario $a$ if in the restricted operation range in the LTGSPLEC model (m$^3$/s). \\
$\underline{u}_{hta\underline{z}}$ & identifier of if the initial volume for hydro $h$ in period $t$ and scenario $a$ is in range $\underline{z}$ or not.\\
$\overline{u}_{hta\overline{z}}$ & identifier of  if the initial volume for hydro $h$ in period $t$ and scenario $a$ is in range $\overline{z}$ or not.\\
$u_{hta}$ & identifier  of if the initial volume for hydro $h$ in period $t$ and scenario $a$ is in the restricted operation range or not.\\

\end{longtable*}

\textbf{Parameters}
\begin{longtable*}{@{}p{1.0cm}p{11.5cm}@{}}
$\mbf{C}_g$ & unitary variable cost of thermal plant $g$ (R\$/MW). \\
$\mbf{CD}$ & cost of deficit (R\$/MW). \\
$\mbf{K}$ & constant that converts water flow (m$^3$/s) to volume (hm$^3$) in a one-month period. \\
$\mbf{NG}$ & number of thermal plants. \\
$\mbf{NH}$ & number of hydro plants. \\
$\mbf{NT}$ & number of months in the planning horizon. \\
$\mbf{NA}_t$ & number of scenarios in period $t$. \\
$\mbf{NZ}_{h}^{\min}$ & number of ranges of volume to determine the minimum outflow of hydro $h$. \\
$\mbf{NZ}_{ht}^{\max}$ & number of ranges of volume to determine the maximum outflow of hydro $h$ in period $t$. \\
$\mbf{D}_{t}$ & demand of active power in period $t$ (MW). \\
$\mbf{Y}_{hta}$ & incremental inflow of hydro $h$ in period $t$ (m$^3$/s) and scenario $a$. \\
$\mbf{\rho}_h$ & productivity constant of hydro $h$ (MWmonth/(m$^3$/s)) \\
$\overline{\mbf{V}}_h$ & maximum storage capacity of the reservoir $h$\\
$\overline{\mbf{ph}}_h$ & maximum generation capacity  for  hydro plants $h$\\
$\overline{\mbf{pg}}_g$ & maximum generation capacity  for  thermal plants $g$.\\
$V^{\rm Res}_{h}$ & percentage of hydro $h$ volume  that defines its  restricted operation range.\\
$V_{h\underline{z}}^{\min}$ & storage level of hydro $h$  for the $\underline{z}$ range of the minimum outflow constraints.\\
$Q^{\min}_{h\underline{z}}$ & minimum outflow of hydro $h$ if the initial storage level  is in $[V_{h\underline{z}}^{\min},V_{h\underline{z}+1}^{\min})$.\\
$V_{ht\overline{z}}^{\max}$ & storage level of hydro $h$ at period $t$ for the $\overline{z}$ range of the maximum outflow constraints.\\
$Q^{\max}_{ht\overline{z}}$ & maximum outflow of hydro $h$ at period $t$ if the initial storage level  is in $[V_{ht\overline{z}}^{\max},
V_{ht \overline{z}+1}^{\max})$.\\
$Q^{\sup} $ & large value compared with possible  outflow values.\\
$T$ & number of seconds per  period.\\
$\underline{a}_{h\underline{z}}$ & angular coefficient of the $\underline {z}$-th part of the function that describes the maximum outflow of the hydro $h$  in the piecewise linear model.\\ 
$\underline{b}_{h\underline{z}}$ & linear coefficient of the $\underline {z}$-th part of the function that describes the maximum outflow of the hydro $h$   in the piecewise linear model.\\
$\overline{a}_{ht\overline{z}}$ & angular coefficient of the $\overline {z}$-th part of the function that describes the maximum outflow of the hydro $h$ in the period $t$ in the piecewise linear model.\\ 
$\overline{b}_{ht\overline{z}}$ & linear coefficient of the $\overline {z}$-th part of the function that describes the maximum outflow of the hydro $h$ in the period $t$  in the piecewise linear model.
\end{longtable*}

\textbf{Sets}
\begin{longtable*}{@{}p{1.0cm}p{11.5cm}@{}}
$\mathcal{G}$ & thermal plants, i.e., $\{1,2,\dots,\mbf{NG} \}$. \\
$\mathcal{H}$ & hydro plants, i.e., $\{1,2,\dots,\mbf{NH} \}$. \\
$\mathcal{T}$ & periods, i.e., $\{1,2,\dots,\mbf{NT} \}$. \\
$\mathcal{M}_h$ & hydro plants upstream of hydro $h$.\\
$\mathcal{A}_t$ & scenarios in  period $t$.\\
$\mathcal{E}$ & hydro plants that must satisfies environmental constraints.\\
$\mathcal{Z}_{h}^{\min}$ & ranges used to define the minimum outflow constraints of hydro $h$, i.e., $\{1,2,\dots, \mbf{NZ}_{h}^{\min}\}$.\\
$\mathcal{Z}_{ht}^{\max}$ & ranges used to define the maximum outflow constraints of hydro $h$ in period $t$, i.e., $\{1,2,\dots, \mbf{NZ}_{ht}^{\max}\}$.\\

\end{longtable*}

\section{Introduction and context}

In Brazil, the hydrothermal operation planning is performed for a time-span of five years, using the {\sc newave} computational model \cite{MacPennaDiniz18} developed by Cepel. One
of the main results of {\sc newave} are the so-called water ``future cost functions'', which represent, by means of piecewise linear models, the monthly opportunity cost of the water stored in reservoirs. The computation of such functions depends on the convexity of the optimization problem solved by {\sc newave} and are a natural output of the application of the well known stochastic dual dynamic programming (SDDP) decomposition strategy \cite{pp91}.

The future cost functions quantify the cost of using water to generate electricity and are used by shorter-term models \cite{DinizCostaMac18}, therefore impacting the entire chain of decisions in the energy business. Inaccurate estimates, i.e., which indicate lower/higher present values of water than the (in practice, unknown) accurate ones leads to an over (under) depletion of reservoirs and, as a consequence, in higher costs of energy production (in the case of cost-minimization problems).

The presence of large river basins stimulates several economic and human activities that may have conflicting objectives among each other. For example, in Brazil, the São Francisco river, which has an extension of over 2,830km, serves as the basis for the development of an entire region of the country. The federal law 9433/97, known as the ``Law of Waters'', promotes a  rational and integrated use of the Brazilian water resources, taking into account multiple uses, such as energy generation, irrigation, navigation, flood control, leisure and tourism, water quality, preservation of aquatic flora and fauna. By this law, in situations of scarcity, water resources should ensure the continuity of supply for essential human activities. The reservoir volume in this case is primarily reserved for human consumption and for watering livestock.

\begin{figure}[hbt]
\centering
\includegraphics{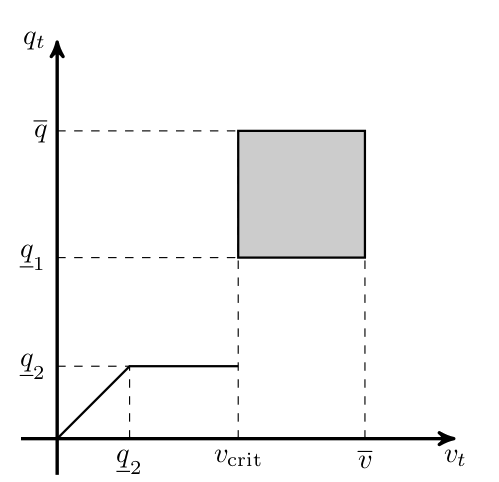}
\caption{Relation between stored volume ($v$) and turbined outflow ($q$) when  ``if-then-else''-like environmental constraints impose a minimum discharge of water.
(source \cite[Fig.5.5]{fgmat})\label{fig1}}
\end{figure}

Since hydroelectric energy represents one alternative use of water,
hydro-generation  must be optimized in conjunction with the other purposes mentioned previously. The application of ``if-then-else'' rules in order to prioritize certain uses of water in the mentioned context may lead to a non-convex and non-connected feasible set for the optimization problem solved by {\sc newave}, as illustrated by Figure~\ref{fig1}.

From the point of view of the optimization problem, this type of constraints can be represented by using some approximations such as piecewise linearizations, or by inserting 0-1 variables in the problem formulation. These modeling variants have different impact on the accuracy of representing the original problem and in the solving techniques employed to calculate the future cost of water. The proposal of the workshop is to approach the theme in two lines of work, by using  {\sc sddp.jl} (\cite{dowson_sddp.jl}, see also \cite{dowson_policy_graph}) and benchmarking their results with the ones obtained by the Newave model \cite{MacPennaDiniz18}. 
On one hand, Julia (together with JuMP \cite{JuMP.jl-2017}) is a package that allows easy modeling of large multi-stage stochastic programming problems and to apply either stochastic dual dynamic programming \cite{pp91} or stochastic integer dual dynamic programming \cite{sddip} to solve the problem. On the other hand, Newave model has been officially used for over 20 years \cite{MacTerry2002} for the official energy planning and price setting of the very large Brazilian system, through application of the SDDP strategy with several improvements along time \cite{MacDuartePenna08,MacPennaDiniz18}. In this sense, during the workshop, comparative studies were carried on, analyzing advantages and disadvantages of different mathematical formulations, particularly with regard to its impact on energy prices and the value of water.

This work is organized as follows. Section~\ref{sec-intro} puts the considered problem in context. Section~\ref{sec-lit} reviews how different hydro constraints have been modelled in the literature, focusing on their impact in operation cost and water value, particularly for long term planning. This section also fixes notation as well as the initial mathematical formulation of dynamic programming approaches that are popular for the problem under consideration. The detailed analysis of an instance of the multi-stage stochastic programming problem is the topic of Section~\ref{sec-ill}. 
By resorting to a simple, yet illustrative, example, several features of the
problem are given gradually, first considering the impact of  minimum/maximum outflow constraints,
and then analyzing how the state-dependence of those constraints makes
the water values nonconvex. Section~\ref{sec-newave} states the mathematical
formulation of the planning problem solved by {\sc newave}.
The different models for environmental constraints that were created  
in the workshop are detailed in Section~\ref{sec-env}. 
The numerical assessment of each considered variant is
presented in Section~\ref{sec-num} and concluding remarks are given in the final
Section~\ref{sec-conc}. Data for the S\~ao Francisco river hydro-bassin
can be found in the Appendix~\ref{sec:test-system-data}.
\section{Generation planning and environmental constraints}\label{sec-intro}

Generation planning/scheduling problems aim to define an operation policy and generation targets for hydro and thermal plants - together with other resources such as wind and solar plants - in order to meet the power demand, while satisfying technical characteristics and constraints of the system components, including the electrical network, on a cost-minimization or a price-maximization framework \cite{WoodWollember96,daSilva01}. Due to the complexity of some of the system constraints, the uncertainty related to natural resources (for example, intermittency of new renewable generation and uncertainty on hydro inflows) and  the large size of real systems, the problem is usually decomposed in long, mid and short term problems, specially in predominantly hydro systems as in Brazil and Norway \cite{Helsethmelo20}.

Among the set of operation constraints that need to be satisfied, we are particularly interested in studying what we label in this document as ``environmental  constraints'' (EC), which are defined as hydraulic constraints for the operation of hydropower and river courses. Such constraints are usually related storage levels of the reservoirs or release (turbined/spilled outflow) of the hydro plants, in order to meet some environmental requirements related to nature, human activities (such as fishing and recreation), or other uses of water, as for example priority water intakes for human/animal, or minimum levels in the watercoures for navigation purposes. The limits established by such constraints may vary along time and may also depend on the current system state, as for example the level of the reservoirs and the hydrological condition of the system.

The modelling of ECs in hydropower planning problems has been of growing interest in the literature in the last decade, due to increasing environmental concerns, and such constraints have become stricter in recent years. In Europe, some revisions have been proposed to hydropower concessions conditions and the implementation of EU Water Framework Directive \cite{Helseth19}. In the case of Brasil, several constraints are enforced by the Water Regulation Agency (ANA) to be considered in the operation of the reservoirs \cite{ONS21}. Some works that describe political/ EC to the watercourse and that may affect optimization of hydropower are discussed in \cite{Helseth19}.

These hydraulic environmental requirements restrict the operation of the hydro plants, specially with the growing share of intermittent renewable generation such as wind or solar, which requires fast ramping resources - such as hydropower - to back up generation in hours with low wind or insollation. More simple EC such as minimum level and maximum ramping of outflow are easier to handle, whereas constraints that impose nonconvex rules for application of constraints that may depend on state variables of the problem (typically, through tables that define minimum/ maximum outflows along intervals for reservoir level) are more involving.

Usually, short term hydrothermal or unit commitment problems are already solved with nonconvex formulations, leading to nonlinear, mixed-integer or mixed-integer linear programming solving strategies \cite{TaktakDambrosi17,ZhengWang15}. In such context, the inclusion of EC is just one more piece that fits the properties of this already complex puzzle, whose main objective is to define generation targets for the day ahead scheduling. On the other hand, in  mid/long term planning problem, the main objective is to define operation policies (in the form of water values) that will guide short term scheduling problems. Such policies are mathematically described as (usually, piecewise linear) cost-to-go functions that are obtained as an output of the application of stochastic programming strategies, such as stochastic dual dynamic programming (SDDP).  In this other context, which is the main interest of this study, convex approximations -- or at least computationally tractable nonconvex formulations -- are desired. We present next a literature review on this subject.

\section{Literature Survey on Environmental Constraints (EC)}\label{sec-lit}

Several hydraulic constraints for power generation planning and scheduling have been proposed in the literature, either properly labeled as  EC  \cite{EdwFlaimHowitt99}, or in a broader context. A far from exhausting list of hydraulic constraints found in the literature is presented below:

\begin{itemize}
    \item time-varying constraints on minimum or maximum release of reservoirs, for each time step; 
    \item time-varying constraints on minimum or maximum storage levels in the reservoirs, for each time step; 
    \item constraints on maximum average water release within a given time window (e.g., a month) \cite{EdwFlaimHowitt99};
    \item constraints on maximum variation of turbined outflow \cite{BorgDambrosio08} or reservoir storages in consecutive time steps \cite{ChangAgan01};
    \item maximum variation in the water level on some river sections downstream to reservoirs \cite{DinizSousa14};
    
\end{itemize}

A more detailed description on those different types of constraints can be seen in \cite{SchafferAdeva20}. We focus our analysis in more recent works whose objective is to assess the modelling challenges as well as the impact - in operation cost and water values - of the inclusion of hydraulic constraints, with emphasis in mid/long term planning problems and in particular, with the application of dynamic programming-like solving strategies.

\subsection{Economic Impact} \label{sec-EconomImpacts-EC}

Reports on the economic impacts of the representation of  EC  in hydropower optimization models for different systems have been surveyed in \cite{Guisandez13}. It addresses minimum outflow and maximum ramping environment constraints in a hydro power optimization problem on a profit maximization framework. The paper applies traditional stochastic dynamic programming with 5 scenarios for each weekly stage, where the subproblems for each stage consists in a MILP problem with hourly discretization. The main objective of the paper is to numerically evaluate the economic losses of the plant owner due to such constraints, with the purpose of supporting a possible  discussion of plant owners with river authorities in order to trade-off the benefits of these constraints with their negative impacts on the profits of the plant owner. However, the paper does not evaluate the impact of environment on the water values themselves, but only on the value of the objective function. Their results (e.g., Fig. 5 in \cite{Guisandez13}) suggest that the higher the level of these constraints, the higher is the derivative of such cost curves. In another paper \cite{Guisandez16a}, the same authors derive regression models for the percentage loss in gain of the hydro producer as a function of the parameters of such constraints. 

\subsection{Impact on water values} \label{sec-WaterValues-EC}

An interesting assessment on the impact of such constraints on the water values themselves is presented in \cite{Guisandez16b}, which derives a mathematical analysis based on the multipliers of the problem constraints. It shows that the water values of the reservoirs increase as the minimum outflow requirement becomes more restrictive (i.e., reaches higher values) and decrease as the maximum ramping rate (in ($m^3/s)/h$) becomes more restrictive (i.e., reaches lower values)). This has an intuitive understanding, since higher minimum releases tend to decrease reservoir levels (because it forces hydro generation or spillage), whereas lower maximum ramping rates limits the flexibility of generation, thus reducing the availability of water. The authors claim further in \cite{Guisandez20} the importance of taking into account such environment constraints in the water value computations for a mid-term planning model to be used in short term simulations. They also show that, for a fixed power capacity of a plant, the impact of these constraints in the water value decreases with the number of hydro units, arguing that the impact for each unit becomes smaller, yielding more flexibility for the hydro plant to handle with such constraints.

\subsection{Modeling of EC in SDP} \label{sec-EC-SDP}

In this section we  briefly describe 
stochastic dynamic programming method (SDP) \cite{BellDrey62} and point out some successful uses of this technique in power generation planning problems.
In order to do this, notice that variables can be separated into two types, \emph{control} and \emph{state} variables, the latter being the result of decisions taken for the
former. Specifically in the application of interest, the reservoir volumes are state variables, while thermal
and hydro generation, turbined and spilt outflows, and deficit are control
variables. This division is used explicitly in SDP-like approaches. In general, consider the abstract formulation presented as follows, where, for simplicity of exposition, we omit the indices related to scenarios along the time steps:

\begin{equation}\label{pbm-abs1}
\begin{array}{clll}
\displaystyle{\min_{x_{[1:T]}, p_{[1:T]}}}&\displaystyle{\sum_{t=1}^{T}}
\left<c_t, p_t\right> &~~&\\
\text{s.t.}&
x_{t} = A_{t-1} x_{t-1}  + B_t p_t + C_t \tilde\xi_t  + d^=_{t}&\mbox{for
}t=1,\ldots,T,& \\
&E_t x_{t} + F_{t} p_t \geq d^\geq_{t}&\mbox{for }t=1,\ldots,T,&
\end{array}
\end{equation}
where
\begin{itemize}
\item[--] $T$ is the number of time steps, possibly large;
\item[--] $\tilde\xi_t$ is a particular realization at time step $t$  of an $M$-dimensional random process. Each process
component $\xi_t(m)$ follows a generalized autoregressive model with time varying order. The realization $\tilde\xi_t$ becomes known at the beginning of time step $t$;
\item[--] $x_t \in X_t\subset \mathbb{R}^{N_x}$ is the state of the system at the end of time step $t$,
with dynamics given by the transition equation and known $x_0$;
\item[--] $p_t \in P_t\subset \mathbb{R}^{N_p}$ is the control variable, applied to the system at time step $t$; and
\item[--] $\left<c_t, p_t\right>$ is the immediate (linear) cost at time step $t$.
\end{itemize}
The dimensions of matrices $A_ t$, $B_t$, $C_t$, $E_t$, $F_t$ 
are, respectively, $N_x\times N_x$, $N_x\times N_p$, $N_x\times M$, $q_t\times N_x$, 
and $N_p\times M$. 
Vectors $d^=_t$ and $d^\geq_t$ are $N_x$- and $q_t$-dimensional. 

By suitably choosing the involved matrices and right hand side terms, 
the abstract inequality constraint includes the equality expressing demand
satisfaction, such as those appearing in energy models, as well as all box constraints.

Since in \eqref{pbm-abs1} the objective and the inequality constraints are
separable by time steps, a sound solving strategy is to apply a stage-wise
decomposition, that we briefly recall below.
More precisely, \eqref{pbm-abs1} is formally equivalent to the following dynamic
programming recursion: having $\mathcal{Q}_{T+1}\equiv 0$,
 solve, for $t=T, T-1,\ldots$:
\begin{equation}\label{policy}
\mathcal{Q}_t(x_{t-1},{\tilde \xi}_{[t-1]})= 
\left\{
\begin{array}{l}
\displaystyle{\min_{x_{t}, p_{t}}} \; \left<c_t,p_t\right> + \mathcal{Q}_{t+1}(x_{t},{\tilde \xi}_{[t]})\\
x_{t} = A_{t-1} x_{t-1}({\tilde \xi}_{[t-1]})  + B_t p_t + C_t \tilde \xi_t  + d^=_{t}\\
E_t x_{t} + F_{t} p_t \geq d^\geq_{t}.
\end{array}\right.
\end{equation}
For these problems, $x_{t-1}({\tilde \xi}_{[t-1]})$ is an optimal value of state $x_{t-1}$ obtained solving a  problem of  form
 \eqref{policy} for time step $t-1$ given the trajectory 
 ${\tilde \xi}_{[t-1]}$ of process $(\xi_{t})$ up to stage $t-1$.
 The first state $x_0$ is known.

 Because setting $\mathcal Q_{T+1}\equiv 0$ makes the
optimization process use all valuable resources (i.e., emptying the reservoirs) at time $T$, the model is said to
suffer from the well-known ``end of the world'' effect. For this reason, the planning
horizon is usually doubled, considering $2T$, and only the output obtained with the
approach until time $T$ is deemed meaningful.

The equivalence between \eqref{policy} and \eqref{pbm-abs1} is formal because
the cost-to-go functions $\mathcal{Q}_{t+1}$ are not known. In the SDP method the state variables are discretized and enumerated so that the future cost function can be exhaustively calculated.

If one applies the SDP to solve the long term planning problem planning, nonconvex constraints can be naturally  included in the problem, since such method does not require the problem to be convex. However, the subproblems themselves derived from the decomposition strategy are more difficult to be solved for being nonconvex, and a discretization of the state space is required, which somehow prevents the application of this method for very large problems and/or with many state variables. 

Regarding the specific subject of this paper, the works \cite{Guisandez16b}, \cite{Guisandez20} have applied the SDP method to consider minimum values and maximum variation of outflow constraints to the reservoirs. Even though the considered problem was  nonconvex, the  EC  were convex and did not depend on the level of the reservoirs. State-dependent environmental constraints in the context of the SDP methodology were included in \cite{SchafferHels22} by considering season-varying volume-dependent maximum inflow constraints for the reservoirs for a mid term planning problem. It considers  quite specific environmental constraints for reservoirs in Norway: during a specific time-window, outflow should be null, unless hydro inflows reach a given threshold, when it becomes limited to a given value, until the reservoir reaches a minimum desired level, when a minimum level constraint takes place. In a different time window, outflow is allowed once the reservoir level does not decrease from one stage to another. Therefore, this is the case of state-dependent constraints both on inflow and reservoir level, which are dealt using the SDP algorithm as an additional state variable that represents the hydrological state of the system. That paper also makes an assessment on the impact of such constraints for the
water values in the reservoirs, showing that, due to the fact that release in
the reservoirs are allowed only when the reservoir level reaches a given
threshold, the water values in the ``with environmental constraint'' are higher
than water values for the ``without environmental constraint'' case, and the cost to go function achieves a non-concave shape (for a profit-maximization problem). The study case is very small (2 reservoirs) and it is applied also using the  JuMP package.

\subsection{Modeling of EC in SDDP} \label{sec-Nonconvex-EC-SDDP}

In order to avoid the huge computational effort to construct future cost functions, in SDDP \cite{pp91} these functions are iteratively approximated in a cyclic process of backward-forward passes. A Benders-like approach is used to replace the (unknown) recourse function
$\mathcal{Q}_{t+1}(x_t,{\tilde \xi}_{[t]})$ with a (known) piecewise-linear
approximation $\mathfrak{Q}_{t+1}(x_t,{\tilde \xi}_{[t]})$,
that is improved along iterations. 
The approximation is the maximum of different cuts, generated by passing
from $t=1$ to $t=T$
different states $(x_{t},{\tilde \xi}_{[t]})$ at each iteration.
Since for linear
programs like \eqref{pbm-abs1}, the cost-to-go functions $\mathcal Q_{t+1}$ 
are indeed piecewise linear, the iterative process eventually terminates.
The final approximation built for the first stage, knowing ${\tilde \xi}_{[1]}$ 
and $x_0$, constitute the FCF that values the opportunity cost of reservoir volume $x_0$.

The  described above is the basis of SDDP.
To deal with a large number of inflow realizations, that yield astronomically large 
scenario trees, SDDP samples trajectories
 ${\tilde \xi}_{[T]}$  of the process $(\xi_{t})$ 
 that are randomly visited at each iteration.
The termination criterion in this case is statistical, based on confidence
intervals, see \cite{shapiro-ejor-2011} for details. A combination of SDP with SDDP has been applied for a long time in Norway \cite{GjelsBels99}.

The two main challenges of including  EC  in SDDP-like solving strategies are discussed in the sequel.

\subsubsection{state dependency}

The first one is ``state-dependency'', which appears when the level of these constraints depend on the values of state variables of the problem, such as reservoir volumes at the beginning of each stage. In this case, the derivatives of such constraints with respect to the value of the state variables should be taken into account when building the Benders cuts of the future cost function in the backward passes of the algorithm. Although it adds some complexity in solving the problems, this is not a major issue, since state dependency in the right hand side of constraints have already been considered for a long time in the SDDP literature \cite{Maceira93,InfangerMorton96,RebenFlack12};

However, if the constraints define a non-convex feasible region, a convexification procedure has to be performed in order to ensure that the resulting future cost function (FCF) remains convex (in the case of cost-minimization) or concave (in the case of profit maximization). For example, the works \cite{Helseth19,SchafferHels22} address state-dependent maximum discharge limits for the hydro plants, discussing which properties should be satisfied in order to allow a convex approximation of these constraints to be considered in the SDDP solving strategy. They arrive to the conclusion that such constraints should have a concave shape. This finding is in line to what has been considered in the {\sc newave} model to represent nonlinear relationship related to state variables in the right hand side of the subproblems \cite{CEPELParab11}: if the parameter $r$ (which defines a physical or operational constraint) of the problem is such that an increase in its value provides an increase (decrease) in the objective function of the problem, the function $r(x)$ that defines the value of the parameter with regard to the state variable $(x)$ (storage in the reservoir) should have a concave (convex) shape, otherwise it has to be transformed in a linear expression. This is the case for the functions that relate the evaporation and minimum outflow constraints as a function of storage, which causes an increase in the objective function and have a concave shape, which, therefore, must be approximated by linear constraints. By contrast, concave expressions of maximum hydro generation with storage (as for example in the hydro production function) and convex expressions of losses due to spillage as a function of water inflows (also used in the modeling of equivalent reservoirs \cite{TcheouCabral12,SagasGuigues12}) allows the use of inequalities approximated by piecewise linear expressions, since the optimal solution will naturally lie in the boundary of these functions.

Another type of  EC  studied in \cite{SchafferHels22} consists in inflow-dependent minimum/maximum discharge constraints, which can also be represented as linear expressions of the inflow, which is an input to the model. However, the particularity of this approach is that, under an autorregressive modeling of inflows \cite{MacBez97}, such relation also creates a time dependency, since the inflow (and as a consequence, the right hand side of  the constraints) depends on the past inflows, which are also state variables of the SDDP approach. This is accounted for by considering an additional term in the computation of the derivatives of the cost to go function, in the same way as {\sc newave} model addresses the existence of other values in the right hand side that depend on the past inflows \cite{DinizCruz21}. Finally, constraints that are dependent on both inflow and storage values are also considered in \cite{Helseth19}, where a three dimensonal linear constraint is proposed to represent this relation. Additional specific constraints that are modeled as virtual reservoirs are also considered in that paper.

\subsubsection{Nonconvex constraints}

A nonconvex nature of some environnmental constraints naturally arises when ``if-then-else'' relations are required. Unfortunately, this is the case for many real applications, since legal requirements tend to express water release constraints as tables instead of mathematical expressions. This requires the consideration of non-convex versions of SDP, as discussed in Section \ref{sec-Nonconvex-EC-SDDP}

The presence of nonconvexities in the SDDP algorithm has been approximately addressed in \cite{PereiraPuschel16} to consider volume-dependent maximum intake of waters in the reservoirs, given by tables, and in \cite{WarlHaug08} to linearize startup costs of thermal units. More sophisticated ways to handle such nonconvexities are proposed in \cite{CerLatRam12}, where a decomposition approach is applied to consider non-convex hydro production function, and in \cite{HelsFodMo16} with a combined SDP/SDDP procedure to consider uncertainty on prices based on a Markov Model. 

\subsection{Modeling of EC in SDDiP} \label{sec-EC-SDDiP}

Instead of applying approximations of the nonconvex  EC  in order to apply the
convex SDDP method, one could explicitly represent such constraints with integer
variables. In this case, the condensed abstract formulation of the problem is as follows:

\begin{equation}\label{pbm-abs2}
\begin{array}{clll}
\displaystyle\min_{\mbox{\scalebox{.6}{$\begin{array}{c}
x_{[1:T]}, p_{[1:T]},
u_{[1:T]}\end{array}$}}}&
\displaystyle{\sum_{t=1}^{T}}
\left< c_t ,p_t \right> &~~&\\
\text{s.t.}&
x_{t} = A_{t-1} x_{t-1}  + B_t p_t + C_t \tilde\xi_t  + d^-_{t}&\mbox{for
}t=1,\ldots,T,& \\
&E_t x_{t} + F_{t} p_t +G_t u_t \geq d^\geq_{t}&\mbox{for
}t=1,\ldots,T,& \\
&u_t\mbox{ has 0-1 components.}&&
\end{array}
\end{equation}
In addition to the objects in  \eqref{pbm-abs1},
we have now matrix $G_t$, of suitable dimensions.

When compared to \eqref{pbm-abs1}, the main difficulty in \eqref{pbm-abs2} is in
the presence of \emph{binary variables}. The SDDP approach relies heavily on
Benders cuts, generated by means of multipliers of the equality constraint. 
However, the combinatorial optimization problem \eqref{pbm-abs2} has no dual
variables.

Currently, the best-established approach to solving the problem is to try to apply the Stochastic Dual Dynamic Integer programming (SDDiP)  method, proposed in \cite{sddip}. In SDDiP, one  first reformulates the nodal problem in a particular way that guarantees strong duality for the Lagrangian dual, and then apply the cuts. In order to  provide the convergence theory, the authors consider that all the state variables are binary, or make a binary expansion of them and use Lagrangian cuts. 

Some variants of the application of this method have been proposed, ranging from: (i) the application of the relaxed continuous version of the subproblems when building the Benders cut - which yields poor and possibly loose approximations of the nonconvex future cost function; (ii) the exact formulation, that assures convergence to the global optimum but requires the use of binary expansion over the continuous state variables in order to obtain a state-space comprised of only binary variables. Several variants in between that yield stronger (tighter) valid approximations of the recourse function have  also been proposed, depending on the type of cut that is built: traditional Benders cuts, Lagrangian cuts or 
Strenghtened Benders cuts.  However, as can be seen in the applications so far \cite{Hjelmeland19,HelsethMo20}, this method lacks maturity for application in large and real multistage hydrothermal planning problems.

In particular, \cite{HelsethMo20} assesses the application of the SDDiP method
for maximum discharge constraints given as a stair-case function of the state
variable related to the storage of the reservoirs, which is very similar to the
ones that have been imposed for the Brazilian system. They have implemented in
the JuMP environemnt the SDDiP method with the three types of cuts presented in
\cite{sddip}. The results show an improvement in the profit results of a final simulation with 1000 scenarios when Lagrangian Cuts are used as compared to 
Strenghtened Benders cuts and then traditional Benders cuts, in this strict order. The convergence results show that the optimality gap is no longer improved when using Benders cuts after a given number of iterations, whereas Strenghtened and Lagrangian cuts are able to keep reducing the gap in further iterations, but with the expense of a much larger CPU time. The study presented in that paper is very similar to what has been proposed in this research, but our aim is to give additional contributions in terms of assessing the impact of the application of  EC  in the water values, and how different solving strategies - including the application of the current convex SDDP approach contributes to obtain accurate evaluation of water values.

\section{EC impact on water values, an illustration}\label{sec-ill}

Prior to more complex analysis taking into account stochastic problems with several reservoirs and time steps, we present below an illustrative analysis of the impact of EC in the water values in the reservoirs. These water values are the derivatives of the FCF, which relates the future cost $FC$ with the storage $V$ in the reservoirs. The FCFs are obtained as an output of the application of dual dynamic programming-based optimization strategies to solve the hydrothermal planning problem. 

We consider a very simple deterministic two-stage problem, whose second-stage subproblem  for scenario $a$ depends on the state variable $\mbf{V_i}$ related to the initial storage in the reservoir for the second stage, which corresponds to the final storage of the first stage. This second stage subproblem, whose solution for several discretized values of $\mbf{V_i}$ defines the FCF for the first stage, is defined as follows:

\begin{mini!}|s|<b>
{}{\sum_{g=1}^{4} \mbf{c}_g pt_{g}  \protect\label{eq:example-fobj}}
{\label{eq:example-model}}{}
\addConstraint{\sum_{g=1}^4 pt_{g} + ph = \mbf{D},} \label{eq:example-demand}
\addConstraint{V_{t+1} + ph = V_{t} + \mbf{Y}_a, } 
\addConstraint{0 \leq V \leq \overline{\mbf{V}}; \ 0 \leq ph\leq \overline{\mbf{ph}},} \label{eq:example-hlim-m1}
\addConstraint{0 \leq pt_{g} \leq \overline{\mbf{pt}}_g,\quad \forall g \in \{1,2,3,4\}, } \label{eq:example-tlim-m1}
\end{mini!}
where, for simplicity of exposition and without loss of generality, storage $V$ is represented in the same unit as hydro generation $ph$. In this sense, minimum/maximum release constraints are equivalent to minimum/maximum generation constraints. The problem data are the following:

\begin{itemize}
    \item two scenario inflows: $Y_1 = 100$ (wet) and $Y_2 = 0$ (dry), with equal probability; 
    \item vector of thermal generation unitary costs: $c = [100,200,300,400]$;
    \item maximum generation capacity: $\overline{\mbf{pt}}_g = 100$ for all thermal plants $g$; 
    \item maximum storage $\overline{\mbf{V}}=400$ and maximum physical discharge $\overline{\mbf{ph}} = 400$ for the reservoir;    
    \item the demand for the second stage is 400;
\end{itemize}

In order to satisfy the relatively complete recourse requirement for the problem (otherwise it would not be possible to evaluate the FCF for some values of storage), we included an unitary penalty value of \$50 for violation of the imposed constraints.

\subsection{Minimum outflow constraints}

\begin{figure}  [hbt]
\centering
\includegraphics[width=1.0\linewidth]{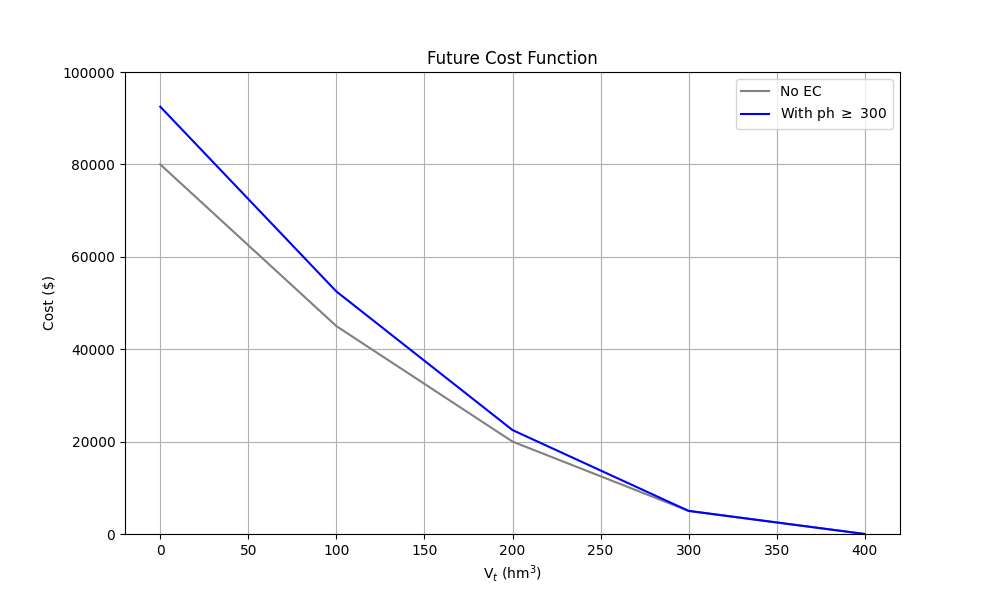}
\caption{FCF for the base case without EC (grey) and with constraints on minimum outflow of the hydro plant (blue)} \label{fig-conceptual-ghmin}
\end{figure}

First we imposed a minimum outflow (generation) constraint $ph \geq 300$, yielding the FCF shown in blue in Figure \ref{fig-conceptual-ghmin}. As expected, future costs are higher with such constraints, for any value of storage. We note that minimum generation constraints lead to higher water values (derivatives in the blue function) as compared to the base case (derivatives of the grey curve). Moreover, this impact occurs specially for lower values of storage, where such constraints cause a higher impact. 

\subsection{Maximum outflow constraints}

Next we imposed a maximum outflow (generation) constraint $ph \leq 250$, yielding the FCF shown in orange in Figure \ref{fig-conceptual-ghmin}. Although costs are also higher with these constraints, we note that the water values are always equal or \textbf{lower} when such constraints is imposed. This effect is more pronounced for higher values of storage, where these constraints are more restrictive, because  prevent the reservoir to take advantage, in the future, of higher values of storage at the end of the first stage. In particular, we can see that the water value is null for values of storage greater than 250, which is the maximum amount of water allowed for generation in the second stage, due to this maximum release constraint.

\begin{figure} [hbt]
\centering
\includegraphics[width=1.0\linewidth]{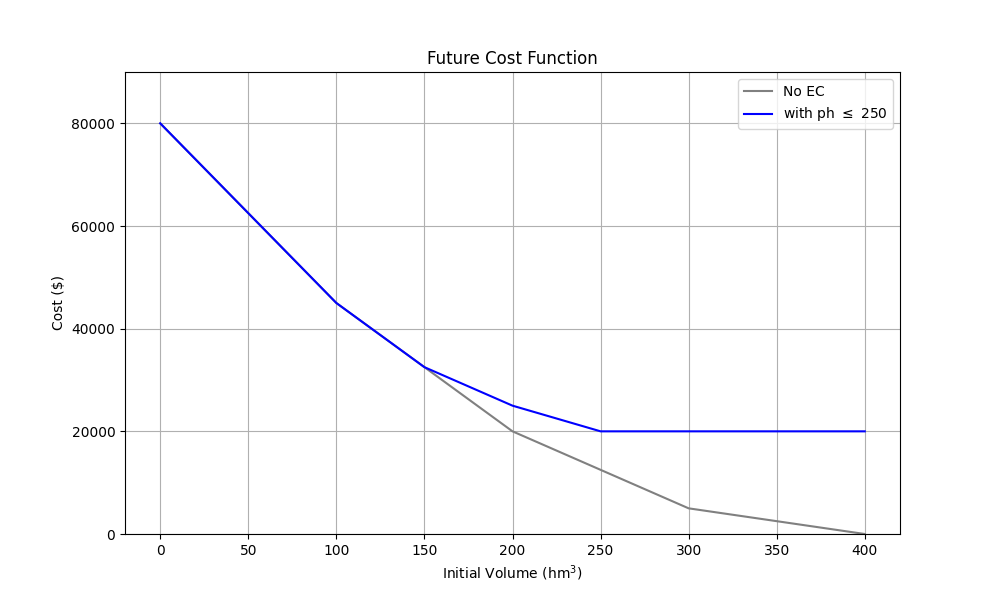}
\caption{FCF for the base case without EC (grey) and with constraints on maximum outflow of the hydro plant (orange)} \label{fig-conceptual-ghmax}
\end{figure}

\subsection{State-dependent Maximum generation constraints}

\begin{figure} [hbt]
\centering
\includegraphics[width=1.0\linewidth]{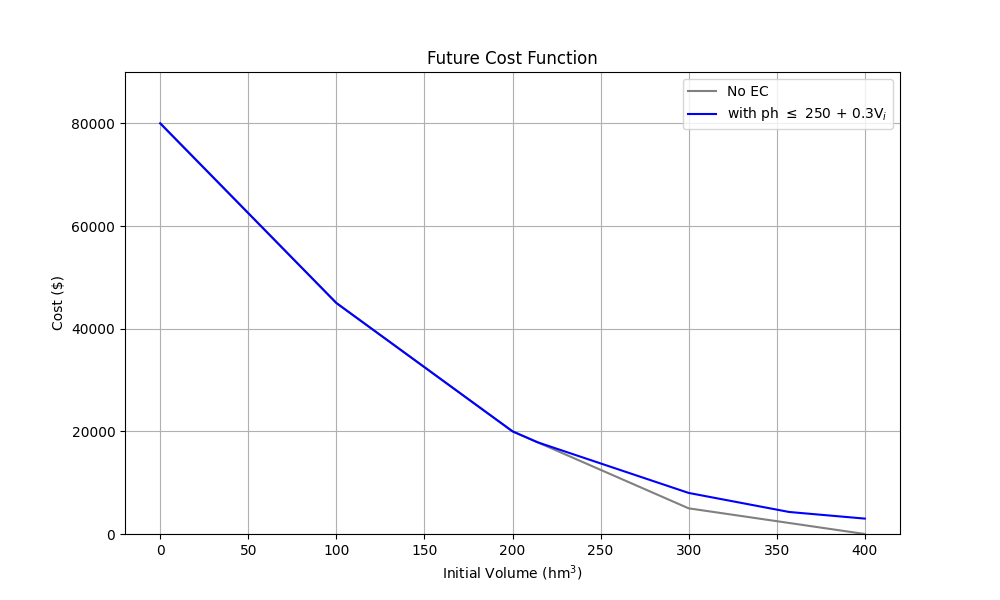}
\caption{FCF for the base case without EC (grey) and with constraints on a state-dependent maximum outflow of the hydro plant (red)} \label{fig-conceptual-ghmax-stateFx}
\end{figure}

We now imposed a state-dependent maximum generation constraint, given by the function $ph \leq 250 + 0.3 V_{i}$, yielding the FCF shown in red in Figure \ref{fig-conceptual-ghmax-stateFx}. Since the maximum outflow increases with storage, such constraints are less restrictive and cause a lower impact in water values as compared to the ones shown in Figure \ref{fig-conceptual-ghmax}. However, these constraints still decrease water values for higher values of storage.

\subsection{Nonconvex Maximum generation constraints}

Finally, we imposed a nonconvex maximum constraint for the second stage, defined as follows:

\begin{itemize}
    \item $ph \leq 250$, for $V_i \leq 300$; 
    \item $ph \leq 10 + 0,8 V_i$, for $V_i > 300$;  
\end{itemize}

\begin{figure} [hbt]
\centering
\includegraphics[width=1.0\linewidth]{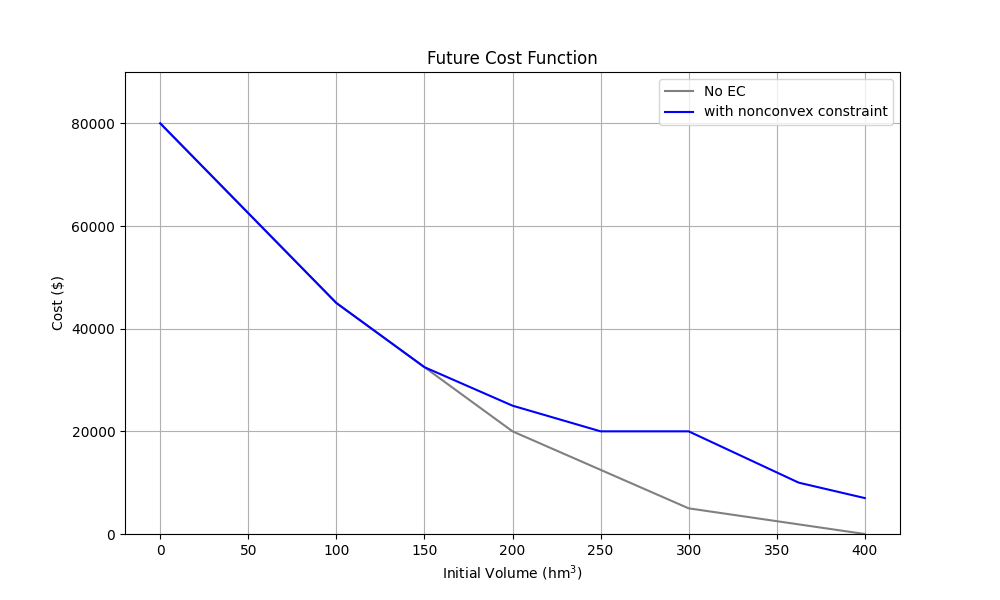}
\caption{FCF for the base case without EC (grey) and with constraints on a mixed fixed / state-dependent maximum generation of the hydro plant (green)} \label{fig-conceptual-ghmax-nonconvex}
\end{figure}

The resulting FCF is presented in green in Figure \ref{fig-conceptual-ghmax-nonconvex}, which shows that those nonconvex constraints cause the FCF to be nonconvex at the region where they have a larger impact, i.e., for higher values of storage.

\section{Mathematical Formulation of the planning problem solved by {\sc
newave}}\label{sec-newave}

The long-term generation scheduling (LTGS) problem considered in this work spans a planning horizon of $\mbf{NT}$ periods, each with a duration of one month, and takes into account the spatial coupling between hydraulic reservoirs. In order to focus on the analysis of this paper, we considered a very simplified version of the problem, where the efficiency of the hydro plants is considered constant and  the demand constraint is modeled as a single bus, i.e., all generations are connected and must meet the demand requirement of a single subsystem. In addition, several modeling features of the program \cite{MacPennaDiniz18} have been disregarded. 

The problem formulation is as follows.

\begin{mini!}|s|<b>
{}{\sum_{t=1}^{\mbf{NT}} \frac{1}{\mbf{NA}_t} \left(\sum_{a \in \mathcal{A}_t}\sum_{g=1}^{\mbf{NG}} \mbf{C}_g pt_{gta}  + \mbf{CD} def_{ta} \right)\protect\label{eq:ltgs-fobj}}
{\label{eq:ltgs-model}}{}
\addConstraint{\sum_{g=1}^\mbf{NG} pt_{gta} + \sum_{h=1}^{\mbf{NH}}ph_{hta} + def_{ta} = \mbf{D}_t, \quad \forall t \in \mathcal{T}, \forall a \in \mathcal{A}_t } \label{eq:ltgs-demand}
\addConstraint{v_{hta} = v_{h,t-1,b(a)} +\mbf{K  Y}_{hta} -\mbf{K}\left(q_{hta} + s_{hta} - \sum_{\forall i\in \mathcal{M}_h}\left(q_{ita}+s_{ita}\right)  \right),}\nonumber
\addConstraint{\quad \forall h \in \mathcal{H}, \forall t \in \mathcal{T}, \forall a \in \mathcal{A}_t  } \label{eq:ltgs-wbal}
\addConstraint{ph_{hta} =  \rho_h q_{hta}, \quad \forall h \in \mathcal{H}, \forall t \in \mathcal{T}, \forall a \in \mathcal{A}_t  } \label{eq:ltgs-ph}
\addConstraint{0 \leq v_{hta} \leq \overline{\mbf{v}_h}; 0 \leq q_{hta} \leq \overline{\mbf{q}_h};\ 0 \leq s_{hta}; \ 0 \leq ph_{hta}\leq \overline{\mbf{ph}_h},} \nonumber
\addConstraint{\quad \forall h \in \mathcal{H}, \forall t \in \mathcal{T}, \forall a \in \mathcal{A}_t } \label{eq:ltgs-hlim}
\addConstraint{0 \leq pt_{gta} \leq \overline{\mbf{pt}}_g,\quad \forall g \in \mathcal{G}, t \in \mathcal{T}, \forall a \in \mathcal{A}_t  } \label{eq:ltgs-tlim}
\addConstraint{def_{ta} \geq 0, \quad  \forall t \in \mathcal{T}, \forall a \in \mathcal{A}_t.  } \label{eq:ltgs-def}
\end{mini!}

The objective function (\ref{eq:ltgs-fobj}) comprises the thermal generation costs over the planning horizon. Equation (\ref{eq:ltgs-demand}) is related to the demand requirement. Equation (\ref{eq:ltgs-wbal}) is the water balance constraint. Equation (\ref{eq:ltgs-ph}) represents the hydro production function, and equations (\ref{eq:ltgs-hlim})-(\ref{eq:ltgs-def}) impose limits on variables, where the notation $\overline{\mbf{x}}$, $\underline{\mbf{x}}$ is the upper and lower limit of variable $x$, respectively

\begin{remark}
\label{rem-abs1}
The LTGS stated in \eqref{eq:ltgs-fobj}-\eqref{eq:ltgs-def} considers one particular instance of the water
inflows $\mbf{K  Y}$ and, for this reason, is a deterministic linear programming problem. 
When inflows are represented by some stochastic process,
the objective function is replaced by an expected cost, as variables and constraints  
now depend on the realizations of the stochastic variable.
The corresponding optimization problem is a multistage stochastic linear
program that we  wrote in  the condensed form \eqref{pbm-abs1}, so that we can apply DP or DDP to solve it (respectively, SDP and SDDP if a stochastic problem is considered).
\end{remark}

\section{ Environmental constraints considered in the model }\label{sec-env}

\begin{figure} [hbt]
\centering
\mbox{\rotatebox{-90}{\includegraphics{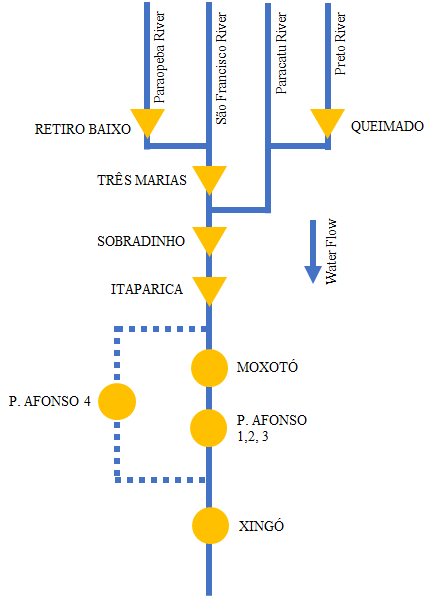}}}
\caption{S\~ao Francisco basin hydro configuration} \label{fig-hydro-conf}
\end{figure}

S\~ao Francisco is one of the most important rivers in Brazil. The
operation of the main reservoirs in its basin is currently carried out by the
Brazilian Independent System Operator (ONS), in accordance with regulations
established by the National Water and Basic Sanitation Agency (ANA).
The schematic diagram in Figure~\ref{fig-hydro-conf} shows the hydro
configuration of the river that is considered in this work.

For some upstream reservoirs, denoted by $h \in \mathcal{E}$, 
ANA specifies three operation modes that control how much water is used to generate
power: restriction, attention, and normal. 
To each mode corresponds an authorized range of operation, depending
on the active storage levels.  The ranges of Tr\^es Marias reservoir
are listed below.
\begin{itemize}

\item Normal operation range: if the active storage level is above
$60\%$ of the total active storage capacity.

\item Attention operation range: if the active storage level is at least
$30\%$ and below $60\%$ of the total active storage capacity.

\item Restricted operation range: if the active storage level is below $30\%$
of the total active storage capacity.
\end{itemize}

For each operation range there is a minimum and maximum total outflow constraint. The minimum total outflow is constant over time, and set to the same value in the 
attention and normal modes.  Therefore, the minimum outflow constraint consists of two
distinct regions, determined by $V^{\rm Res}_h$, which is the threshold level between the attention and restricted ranges. 
Três Marias reservoir has a minimum outflow of 100 $m^3/s$ if the attention operation range is active, and 150 $m^3/s$ otherwise. For each hydro plant $h$, we will denote by $Q^{\min}_{h1}$ the minimum outflow in the restricted operation range and by $Q^{\min}_{h2}$ otherwise. In particular, for Três Marias reservoir, $V^{\rm Res}= 0.3$  and $Q^{\min}=[100,150]$.

\subsection{First EC model}\label{sec:modelA}
In order to model this situation we will use the binary variable $u_{hta}$, for every hydro plant $h$, time $t$ and scenario $a$, to identify if we are in the restriction zone or not. When the reservoir level $h$ is between 0 and $V^{\rm Res}_h$, we must have $u_{hta}=1$, while for the level greater than $V^{Res }_h$ we must impose that $u_{hta}=0$. To avoid using open intervals, let us assume that either of the two possibilities of $u_{hta}$ can be used if the reservoir level is exactly at the $V^{\rm Res}_h$ value. Considering this we have:
 \begin{equation}\label{mindefluselec}
 V^{\rm Res}_{r(h)} (1- u_{hta}) \le \frac{v_{r(h),t-1,b(a)}}{\overline{\mbf{V}}_{r(h)} }\le V^{\rm Res}_{r(h)} u_{hta} + (1 - u_{hta}), \quad \forall h \in \mathcal{E}, \forall t \in \mathcal{T}, \forall a \in \mathcal{A}_t, 
 \end{equation}
  \begin{equation}\label{binmin}
 u_{hta} \in \{0,1\}, \quad \forall h \in \mathcal{E}, \forall t \in \mathcal{T}, \forall a \in \mathcal{A}_t.
\end{equation}
Once it is identified whether we are in the restricted operation range or not, we must impose that the outflow is greater than the specified limit. This can be written as
$$
 Q^{\min}_{h1} u_{hta}+  Q^{\min}_{h2}(1-u_{hta})\leq q_{hta}+s_{hta}, \quad \forall h \in \mathcal{E}, \forall t \in \mathcal{T}, \forall a \in \mathcal{A}_t.
$$
Constraint infeasibility is handled by adding a deficit variable (i.e., a slack variable for constraint violation), as follows:
\begin{equation}\label{mindefluval}
 Q^{\min}_{h1} u_{hta}+  Q^{\min}_{h2}(1-u_{hta})\leq q_{hta}+s_{hta}+def^{\min}_{hta}, \quad \forall h \in \mathcal{E}, \forall t \in \mathcal{T}, \forall a \in \mathcal{A}_t,
\end{equation}
\begin{equation}\label{defminmaior0} 
 def^{\min}_{hta} \geq 0,  \quad \forall h \in \mathcal{E}, \forall t \in \mathcal{T}, \forall a \in \mathcal{A}_t. 
\end{equation}

We now turn our attention to ANA's constraints on maximum outflow.
For the restricted operation range, the total maximum outflow is
determined by ONS following ANA recommendation and, for the normal operation
range, there is no maximum outflow. The maximum outflow within the attention
operation range, however, depends on values established by ONS for each month,
based on the initial active storage level of the
reservoir. Table~\ref{table:tres_marias} shows the maximum outflow to be
considered for the Tr\^es Marias reservoir in each month (from December 2022 to
November 2023), depending on its active storage level at the beginning of that
month.

\begin{table}
\centering
\begin{tabular}{|c|cccccccc|}
\hline
           & 150     & 200     & 250     & 300     & 350     & 400     & 450     & 500\\
           & $m^3/s$ & $m^3/s$ & $m^3/s$ & $m^3/s$ & $m^3/s$ & $m^3/s$ & $m^3/s$ & $m^3/s$\\
\hline
01/12/2022 & --      & 30.0\%  &  33.9\% &  39.1\% &  44.2\% &  49.4\% &  54.5\% &  59.7\%\\
01/01/2023 & 30.0\%  & 33.6\%  &  38.7\% &  43.8\% &  48.9\% &  54.1\% &  59.2\% &  --    \\
01/02/2023 & 32.0\%  & 37.1\%  &  42.3\% &  47.4\% &  52.5\% &  57.7\% &  --     &  --    \\
01/03/2023 & 34.1\%  & 39.3\%  &  44.5\% &  49.8\% &  55.0\% &  --     &  --     &  --    \\
01/04/2023 & 38.1\%  & 43.3\%  &  48.5\% &  53.7\% &  58.9\% &  --     &  --     &  --    \\
01/05/2023 & 40.4\%  & 45.6\%  &  50.8\% &  56.1\% &  --     &  --     &  --     &  --    \\
01/06/2023 & 40.0\%  & 45.2\%  &  50.4\% &  55.6\% &  --     &  --     &  --     &  --    \\
01/07/2023 & 32.5\%  & 37.7\%  &  42.9\% &  48.2\% &  53.4\% &  58.6\% &  --     &  --    \\
01/08/2023 & 30.0\%  & 34.0\%  &  39.2\% &  44.4\% &  49.7\% &  54.9\% &  --     &  --    \\
01/09/2023 & --      & 30.0\%  &  31.3\% &  36.4\% &  41.6\% &  46.8\% &  51.9\% &  57.1\%\\
01/10/2023 & --      & --      &  --     &  30.0\% &  32.5\% &  37.7\% &  42.9\% &  48.1\%\\
01/11/2023 & --      & --      &  30.0\% &  30.7\% &  35.9\% &  41.0\% &  46.2\% &  51.3\%\\
\hline
\end{tabular}
\caption{Maximum outflow for Tr\^es Marias reservoir based on its initial active storage level.}
\label{table:tres_marias} 
\end{table}

In this table, two adjacent storage levels (in the same row) define a
closed-open interval. If the active storage level of the reservoir lies into
that interval, then the maximum outflow is to be taken as the one associated
with the first endpoint of the interval. For instance, for February 2023, if
the active storage level is $45\%$ the total active storage capacity of the
reservoir, then it lies in the interval $[ 42.3\% , 47.4\% )$ and the maximum
outflow to be considered must be $250 m^3/s$.

Maximum outflows are also determined for run-of-the-river power plants based on
the active storage level of some upstream reservoir. For instance, the maximum
outflow of Xing\'o power plant is determined based on the active storage level
of Sobradinho reservoir. Furthermore, the maximum outflow of Sobradinho
reservoir is given by the maximum outflow of Xing\'o power plant plus a
constant (which depends on the month). In general, a reservoir $h$ maximum outflow  is determined based on the active storage level of a reservoir whose index is  $r(h)$ 
(the reference reservoir) at the beginning of the month.

Let $V^{\max}_{h\overline{z}}$,  $\overline{z} \in \mathcal{Z}_{h}^{\max} \equiv \{1,\ldots,\mbf{NZ}_{h}^{\max}+1\}$, be the values that divide the storage levels of reservoir $h$. If
the active storage level of the  reference reservoir of hydro $h$ lies in the interval $[ V_{r(h)\overline{z}}^{\max} ,
 V_{r(h)\overline{z}+1}^{\max})$, then the maximum outflow to be considered for reservoir $h$ must
be $Q^{\max}_{h\overline{z}}$. In order not to deal with open intervals, we simplify and
consider only closed intervals, so that the interval $[ V_{h\overline{z}}^{\max} ,
 V_{h\overline{z}+1}^{\max})$
is replaced by $[ V_{h\overline{z}}^{\max} ,
 V_{h\overline{z}+1}^{\max}]$. Thus, notice that, if the initial
active storage level lies in the frontier of an interval, there could be two
possible ranges for that initial volume.

Consider once more the Tr\^es Marias reservoir. It is its own reference
reservoir (so we would have $r(h) = h$). In the restricted operation range,
when the active volume of the reservoir is below $30\%$, the maximum outflow should be determined by ANA. Since this value is not available, we chose to consider it as  $100 m^3/s$, the minimum outflow enforced by the regulations.  For the normal operation range, where the outflow  is unlimited, we just chose a sufficiently large value $Q^{\sup}$. Thus, for May 2023, we would have for this reservoir that
$$V^{\max}=[0,0.3,0.404,0.456,0.508,0.561,0.6,1]$$   and the outflow limit would be 
$$Q^{\max} =[100, 150,150,200,250,300,Q^{\sup}].$$ 

Now, let us define the binary variable $\overline{u}_{hta\overline{z}}$ to identify if the initial volume of the reference reservoir of hydro $h$,  in period $t$ and under scenario $a$, is in the  $\overline{z}$-th interval $[ V_{ht\overline{z}}^{\max}, V_{ht\overline{z}+1}^{\max}]$. So the constraints that  determine the maximum outflow are:
\begin{equation}\label{maxdefluselec}
    V_{r(h)t\overline{z}}^{\max} \overline{u}_{hta\overline{z}}\le \frac{v_{r(h),t-1,b(a)}}{\overline{\mbf{V}}_{r(h)} } \le V_{r(h)t\overline{z}+1}^{\max} \overline{u}_{hta\overline{z}} + (1 -\overline{u}_{hta\overline{z}} ), \quad \forall h \in \mathcal{E}, \forall t \in \mathcal{T}, \forall a \in \mathcal{A}_t, 
\end{equation}
\begin{equation}\label{binmax}
\overline{u}_{hta\overline{z}} \in \{0,1\}, \quad \forall h \in \mathcal{E}, \forall t \in \mathcal{T}, \forall a \in \mathcal{A}_t, \forall \overline{z} \in \mathcal{Z}_{h}^{\max}. \end{equation}
To enforce that the initial storage level be associated with only one
interval, we must have
\begin{equation}\label{maxdeflu1u}
\sum_{\overline{z} \in \mathcal{Z}_{h}^{\max}} \overline{u}_{hta\overline{z}}  = 1, \quad \forall h \in \mathcal{E}, \forall t \in \mathcal{T}, \forall a \in \mathcal{A}_t.
\end{equation}

Thus, the maximum outflow for reservoir $h$, in period $t$ and under scenario $a$ is given by
$$
 q_{hta}+s_{hta} \leq  \sum_{\overline{z} \in \mathcal{Z}_{h}^{\max}} \overline{u}_{hta\overline{z}} Q^{\max}_{r(h)t\overline{z}}.
$$
In order to ensure the feasibility of maximum outflow, a slack variable is introduced, resulting in the constraints
\begin{equation}\label{maxdefluval}
 q_{hta}+s_{hta} \leq  \sum_{\overline{z} \in \mathcal{Z}_{h}^{\max}} \overline{u}_{hta\overline{z}} Q^{\max}_{r(h)t\overline{z}}+def^{\max}_{hta}, \quad \forall h \in \mathcal{E}, \forall t \in \mathcal{T}, \forall a \in \mathcal{A}_t,
\end{equation}
\begin{equation}\label{defmaxmaior0}
    def^{\max}_{hta} \geq 0, \quad \forall h \in \mathcal{E}, \forall t \in \mathcal{T}, \forall a \in \mathcal{A}_t.  
\end{equation}

Therefore, we define the LTGSEC model as the LTGS model plus  EC . To define its objective function we denote by $\mbf{CD}^{\min}$ and $\mbf{CD}^{\max}$ the cost of the deficit of the minimum and maximum outflow constraints, respectively. So the aim in  the LTGSEC model is
\begin{equation}\label{eq:ltgsec-fobj} 
 \min   \sum_{t=1}^{\mbf{NT}} \frac{1}{\mbf{NA}_t} \left(\sum_{a \in \mathcal{A}_t}\sum_{g=1}^{\mbf{NG}} \mbf{C}_g pt_{gta}  + \mbf{CD} def_{ta} +    \sum_{h \in \mathcal{E}} \left(\mbf{CD}^{\min} def_{hta}^{\min}+\mbf{CD}^{\max} def_{hta}^{\max}   \right) \right)
\end{equation}
subject to the constraints \eqref{eq:ltgs-demand} - \eqref{eq:ltgs-tlim} and \eqref{mindefluselec} - \eqref{defmaxmaior0}.

\begin{remark} \label{rem-abs2}
When compared to the LTGS model, we now have to deal with
the binary variables that discriminate the different operational modes
(denoted by $u$).  Once again, we are interested in the stochastic version of the problem. In this way, the model can be seen as a particular case of the abstract form \eqref{pbm-abs2}, and so, using a binary expansion of the state variables, we can apply SDDiP to solve it.  However, the discretization of state variables leads to a substantial increase in binary variables, impairing computational performance. More than that, in \cite{sddip} the cuts that proved to be most efficient in practice are not the ones that support the convergence theory. In this way, an approach without guaranteed convergence, but which demonstrates good computational performance, may be interesting. This is the case of a heuristic implemented in  the package  {\sc sddp.jl}, where integers variables are relaxed and handled as the continuous ones, and there are some ways to compute subgradients in the backward pass.
We denote the use of this strategy by iSDDP, which is similar to pricing approaches in unit-commitment problems, see \cite{oneill-etal-2005}, \cite{cepel-2020}.

\end{remark}

\subsection{Alternative model - piecewise linear approximations}
\label{sec:modelB}

In this second proposal we make an extension of the previous formulation, allowing the maximum and minimum flow restrictions to be piecewise linear but not necessarily constant in each interval. In this way, we can adjust the slope of each range so that the  constraints of the model  become continuous. Furthermore, in regions where the maximum outflow constraint is concave, it is possible to eliminate the integer variables that determine the partition of this region, since it can be replaced by the requirement that the outflow be smaller than all the lines that describe the maximum outflow of that region.

To build a piecewise linear model, we linearly interpolate the outflow values indicated for the beginning of each interval of the attention range. When analyzing the lines that interpolate these data, we observe that the data are almost linear, except for the first interval. In addition, we verified that the angular coefficient of the lines had an increase of at most $2\%$, which confirms that a concave approximation is satisfactory. With this information in hand, we adjusted the angular coefficients, extending the straight line of the previous interval when we obtained an increase in the slope, so that the model becomes in fact concave.

For the normal operation range, we have no outflow limitation by ANA regulations. However, we believe that the impact would be quite limited if we included the constraints resulting from the last division of the attention operation range for the normal range as well. This is because the constraint is not very tight for this zone and maximum outflow increases proportionally to the volume of the reservoir, so we would have reasonably high limits. On the other hand, in the restriction range the maximum flow is constant. In this way, we would only have the concavity of the constraints if the maximum outflow is non-increasing in the attention operation range, which is unreasonable. Another possibility of trying to make a concave approximation of the restriction would be to propose an affine function to describe the maximum flow also in the restricted operation range. This alternative is not successful too, as any slope that is greater than those consistent with growth in the attention range results in negative maximum outflows for some volume values in the restriction range. Thus, there seems to be no way to avoid using an integer variable to identify whether or not the volume level of a certain reservoir $h$ is less than the value $V^{\rm Res}_h$, which determines the restriction range.

This way, as we did with the minimum outflow constraints in the LTGS model, for every hydro plant $h$, time $t$ and scenario $a$, we have a single binary variable $u_{hta}$ to decide if we are in the restriction zone ($u_{hta}=1$) or not ($u_{hta}=0$). So the 
 constraints  that identify if we are in restricted operation range or not become
 \begin{equation}\label{maxdefluselecPL}
 V^{\rm Res}_{r(h)} (1- u_{hta}) \le \frac{v_{r(h),t-1,b(a)}}{\overline{\mbf{V}}_{r(h)} }\le V^{\rm Res}_{r(h)} u_{hta} + (1 - u_{hta}), \quad \forall h \in \mathcal{E}, \forall t \in \mathcal{T}, \forall a \in \mathcal{A}_t,
\end{equation}
  \begin{equation}\label{binmaxPL}
 u_{hta} \in \{0,1\}, \quad \forall h \in \mathcal{E}, \forall t \in \mathcal{T}, \forall a \in \mathcal{A}_t.
\end{equation}

 If we are not in the restricted operation range, we limit the outflow by the previously discussed affine functions, whose coefficients we will denote by $\overline{a}_{ht\overline{z}}$ and $\overline{b}_{ht\overline{z}}$. Using a sufficiently large $Q^{\sup}$ value, we guarantee that if we are in the restriction range, these constraints are trivially satisfied. On the other hand, in the restriction zone, we limit the outflow by a fixed constant  value $Q^{\max}_{ht1}$. Adding slack variables $def^{\max}_{hta}$ and $def^{\rm Res}_{hta}$, these sets of constraints are described, respectively, as follows
  \begin{equation} \label{maxoutflowatencao}
 q_{hta}+s_{hta} \le \overline{a}_{ht\overline{z}} v_{r(h),t-1,b(a)}+\overline{b}_{ht\overline{z}}+Q^{\sup}u_{hta}+ def^{\max}_{hta}, \quad \forall h \in \mathcal{E}, \forall t \in \mathcal{T}, \forall a \in \mathcal{A}_t,
    \end{equation}
      \begin{equation} \label{defmaxmaior0PL}
def^{\max}_{hta} \geq 0, \quad \forall h \in \mathcal{E}, \forall t \in \mathcal{T}, \forall a \in \mathcal{A}_t
    \end{equation}
 and   
 \begin{equation} \label{maxoutflowretricao}
q_{hta}+s_{hta} \le Q^{\max}_{ht1}+Q^{\sup}(1-u_{hta})+def^{\rm Res}_{hta}, \quad \forall h \in \mathcal{E}, \forall t \in \mathcal{T}, \forall a \in \mathcal{A}_t,
    \end{equation}
 \begin{equation} \label{defResmaior0PL}
def^{\rm Res}_{hta} \geq 0, \quad \forall h \in \mathcal{E}, \forall t \in \mathcal{T}, \forall a \in \mathcal{A}_t.
    \end{equation}

Let us consider now the minimum outflow constraints. Once again, our first intention would be to make a formulation that avoids the use of integer variables. For this, it would be necessary that the function that describes the minimum flow of the reservoir as a function of its volume be convex. However, we must have a non-decreasing function, which is incompatible with the fact that the function is constant for volume values above a certain threshold. In this way, we decided to keep the integer variables that determine the intervals in the definition of the minimum flow constraint. 

However, due to natural physical reasons, this constraint cannot be always
satisfied. When the volume of the reservoir is too low, and the inflow is not
large enough, there is not enough water to satisfy the minimum outflow
requirement. For instance,  consider the case in which the reservoir is empty and the inflow is $60
m^3/s$. If the minimum outflow is $100 m^3/s$, the best one can do is to
discharge the total inflow, i.e., $60 m^3/s$. The $100 m^3/s$ minimum outflow
is therefore a requirement that may not be physically possible to fulfill and
may not be part of the mathematical model. What can be enforced though is
to discharge the maximum possible amount of water when the discharge cannot
reach its minimum requirement. Since the inflow is stochastic, we decided to enforce the minimum outflow to be at least the initial volume of the reservoir. 
The advantage of such modeling choice is that it eliminates slack variables.
This is attractive because penalizing slack variables in the
objective function is somehow artificial. On the other hand, this relaxation may favor situations where the level of the reservoirs is below the restriction level. We show in our numerical experiments
the impact of those penalties in the final output of the model.

In order to model this situation, we again make a partition of the reservoir levels into $\mbf{NZ}^{\min}_h$ intervals, now denoted by $[ V_{h\underline{z}}^{\min}, V_{h\underline{z}+1}^{\min}] $, and define a piecewise linear function on each of them.
Similar to what we done before,  we use binary variables $\underline{u}_{hta\underline{z}}$ to define if the reservoir $h$ level volume in time $t$ and under scenario $a$ is in the $\underline{z}$-th interval. Moreover, we denote by $\underline{a}_{h\underline{z}}$ and $\underline{b}_{h\underline{z}}$ the slopes and the linear coefficients of the affine functions in each part. 

Let $T$ be  the number of seconds per period. Then, the amount of water associated with a volume $v_{h,t-1,b(a)}$ that
can be discharged per second in a given period $t$ is $\frac{v_{h,t-1,b(a)}}{T}$. Thus, the minimum
outflow constraint  should be replaced by the minimum  between this value and the one established in ANA's technical regulations. So, for $v_{h,t-1,b(a)}<Q^{\min}_{h1}T$ we have an adjustment of the constraint imposed by the ANA. Therefore, we have that $V^{\min}_{h1}=0$, $V^{\min}_{h1}=Q^{\min}_{h1}T$, 
$V^{\min}_{h2}$ is the threshold level that defines the restricted operation range and $V^{\min}_{h3}=1$. Moreover, $\underline{a}_{h1}=\frac{1}{T_t}$, $\underline{a}_{h2}=\underline{a}_{h3}=\underline{b}_{h1}=0$, $\underline{b}_{h2}=Q^{\min}_{h1}$ and $\underline{b}_{h3}=Q^{\min}_{h2}$. Although we fixed both the number of intervals and the way of choosing affine functions in our model, we decided to present it in general form because we believe that other choices of piecewise linear functions would also be interesting to be studied in the future.

The following constraints identify which range the initial volume level is located:
\begin{equation}\label{mindefluselecPL}    V_{ht\underline{z}}^{\min} \underline{u}_{hta\underline{z}}\le \frac{v_{h,t-1,b(a)}}{\overline{\mbf{V}}_{h} } \le V_{ht\underline{z}+1}^{\min} \underline{u}_{hta\underline{z}} + (1 -\underline{u}_{hta\underline{z}} ), \quad \forall h \in \mathcal{E}, \forall t \in \mathcal{T}, \forall a \in \mathcal{A}_t, 
\end{equation}
\begin{equation}\label{binminPL}
\underline{u}_{hta\underline{z}} \in \{0,1\}, \quad \forall h \in \mathcal{E}, \forall t \in \mathcal{T}, \forall a \in \mathcal{A}_t, \forall \underline{z} \in \mathcal{Z}_{h}^{\min}. \end{equation}
Next we  enforce that the initial storage level is associated with only one
interval 
\begin{equation}\label{mindeflu1uPL}
\sum_{\underline{z} \in \mathcal{Z}^{\min}} \underline{u}_{hta\underline{z}}  = 1, \quad \forall h \in \mathcal{E}, \forall t \in \mathcal{T}, \forall a \in \mathcal{A}_t.
\end{equation}
Finally, the minimum outflow is defined by
\begin{equation}\label{mindefluvalPL}
Q^{\sup}(\underline{u}_{hta\underline{z}}-1)+\underline{a}_{h\underline{z}} v_{h,t-1,b(a)}+\underline{b}_{h\underline{z}} \le  q_{hta}+s_{hta}, \quad \forall h \in \mathcal{E}, \forall t \in \mathcal{T}, \forall a \in \mathcal{A}_t,
\end{equation}
where the term $Q^{\sup}(\underline{u}_{hta\underline{z}}-1)$ makes the constraint trivially satisfied when $\underline{u}_{hta\underline{z }}=0$.
So, the alternative model of LTGSEC,  which includes  EC  using piecewise linear functions,  is defined by
\begin{equation}\label{eq:ltgsplec-fobjPL} 
 \min   \sum_{t=1}^{\mbf{NT}} \frac{1}{\mbf{NA}_t} \left(\sum_{a \in \mathcal{A}_t}\sum_{g=1}^{\mbf{NG}} \mbf{C}_g pt_{gta}  + \mbf{CD} def_{ta} +  \sum_{h \in \mathcal{E}} \mbf{CD}^{\max} def_{hta}^{\max}    \right)
\end{equation}
subject to the constraints \eqref{eq:ltgs-demand} - \eqref{eq:ltgs-tlim} and \eqref{maxdefluselecPL} - \eqref{mindefluvalPL}.

\begin{remark}
    The stochastic version of this variant still fits the abstract formulation \eqref{pbm-abs2}, and can be solved by SDDiP and by iSDDP.
\end{remark}

\bigskip

In the numerical assessment that follows, the LTGSEC formulations described in Sections \ref{sec:modelA} and \ref{sec:modelB} are referred to as \emph{Variant A}
and \emph{Variant B}, respectively.

\section{Numerical tests}\label{sec-num}

To evaluate the different proposals presented in this report, the experiments were conducted on a test system containing a subset of hydrothermal plants from the northeastern Brazilian subsystem. The data from this system can be found in Section~\ref{sec:test-system-data}.

The experiments were performed in a computer with a Intel i7 12700H CPU, 16 GB of RAM running Windows operational system. This problem was written in Julia using the JuMP package~\cite{JuMP.jl-2017} and solved using the {\sc sddp.jl} package~\cite{dowson_sddp.jl} in Julia with the Gurobi~\cite{gurobi} solver.

This section is organized as follows. In Section~\ref{sec:newave}, we present a comparison of results obtained by of our implementation of the model and solution procedure in Julia with that obtained by {\sc newave}. In Section~\ref{sec:impact_ec}, we analyse the impacts of the environmental constraints considered in this report. The remaining results concerning the Zinf of all cases are presented in Appendix.

\subsection{System model consistency w.r.t. {\sc newave}}
\label{sec:newave}

To ensure the validity of our proposed optimization model and the efficacy of the SDDP problem-solving approach, we implemented the optimization problem without incorporating EC. Our initial implementation considered the following attributes:

\begin{itemize}
    \item Study period: January 2021 to December 2021 (12 months).
    \item Duration of each period: 730 hours.
    \item Number of openings used in the backward pass:
    \begin{itemize}
        \item First month: 1
        \item Other months: 4
    \end{itemize}
    \item Number of series used in the final simulation: 100.
\end{itemize}

The number of series applied during the forward pass of SDDP varied between the {\sc newave} and our Julia implementation. For {\sc newave}, we considered 10 simulations during the forward pass. In contrast, our Julia implementation only involved a single simulation in the forward pass.

\begin{figure} [hbt]
\centering
\includegraphics[width=1.00\textwidth]{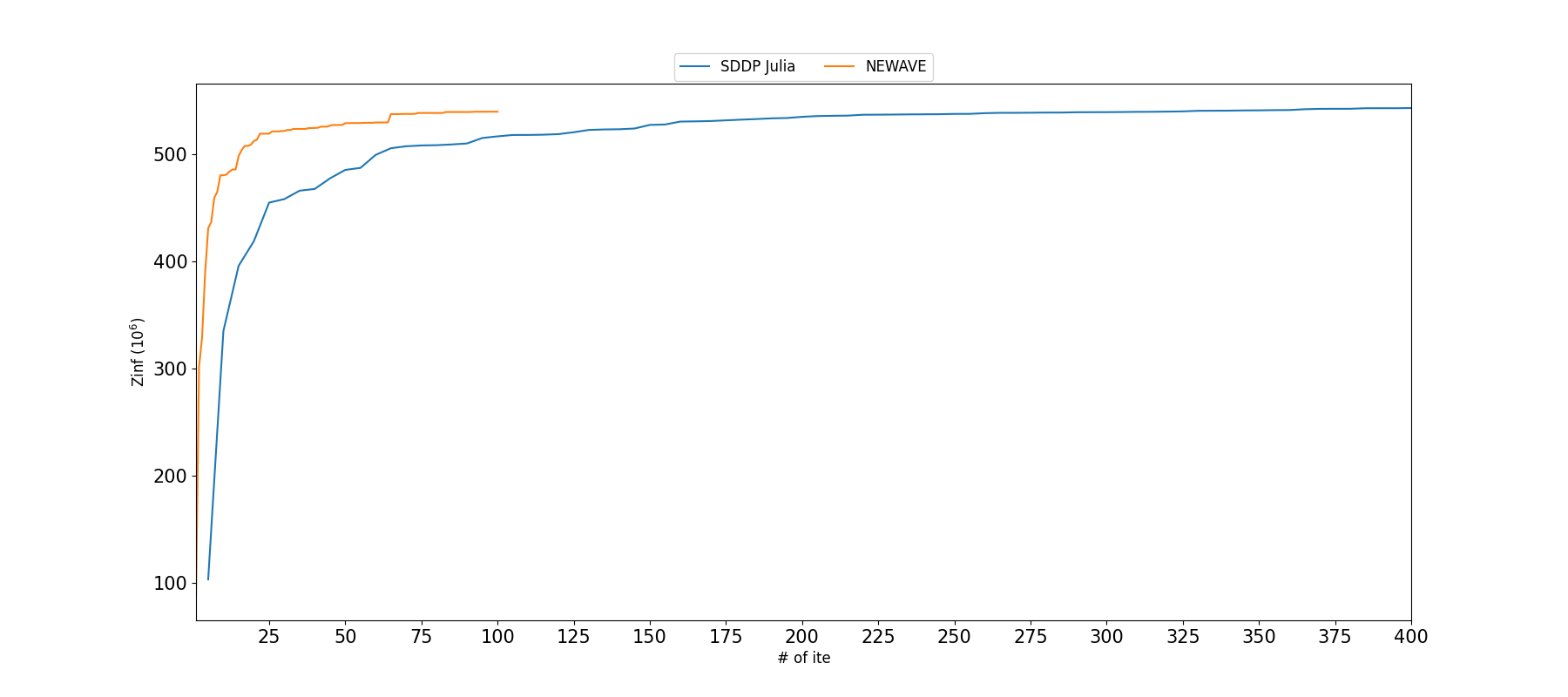}
\caption{Evolution of the lower bound with both computational tools} \label{fig-zinf-case01}
\end{figure}

\begin{figure} [hbt]
\centering
\includegraphics[width=1.00\textwidth]{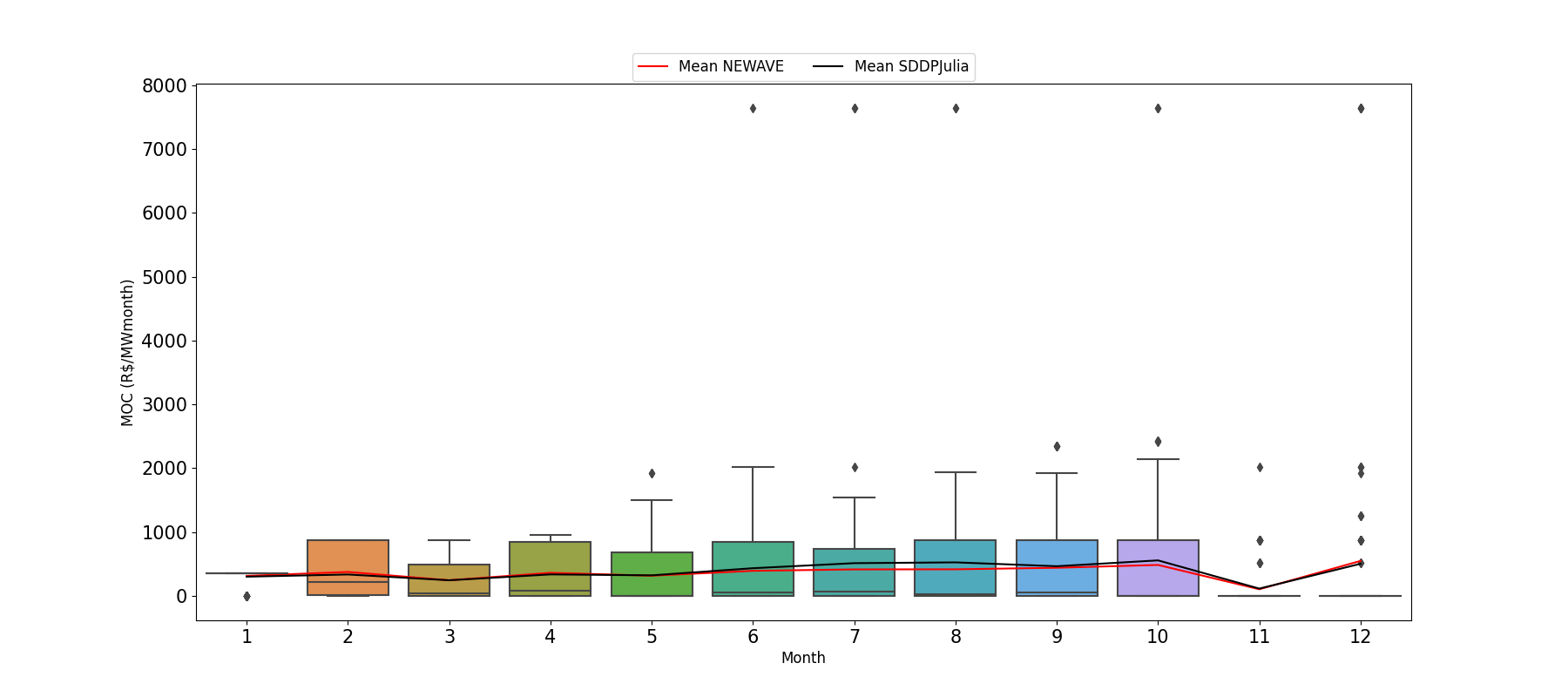}
\includegraphics[width=1.00\textwidth]{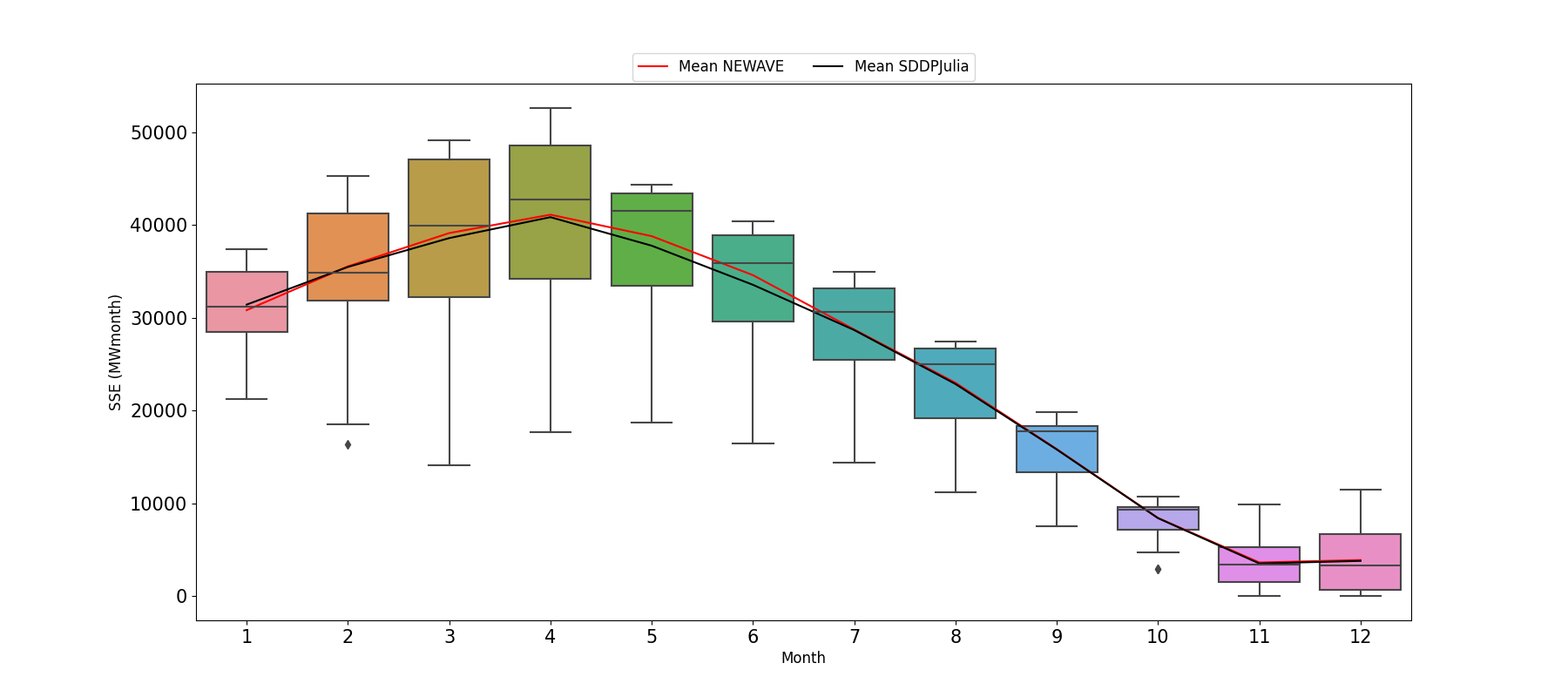}
\caption{Monthly marginal operation cost and system stored energy,
mean and boxplots over the same 100 simulations for {\sc sddp.jl} and {\sc newave}.} 
\label{fig-cmo-ear-case01}
\end{figure}

\begin{figure} [hbt]
\centering
\includegraphics[width=1.00\textwidth]{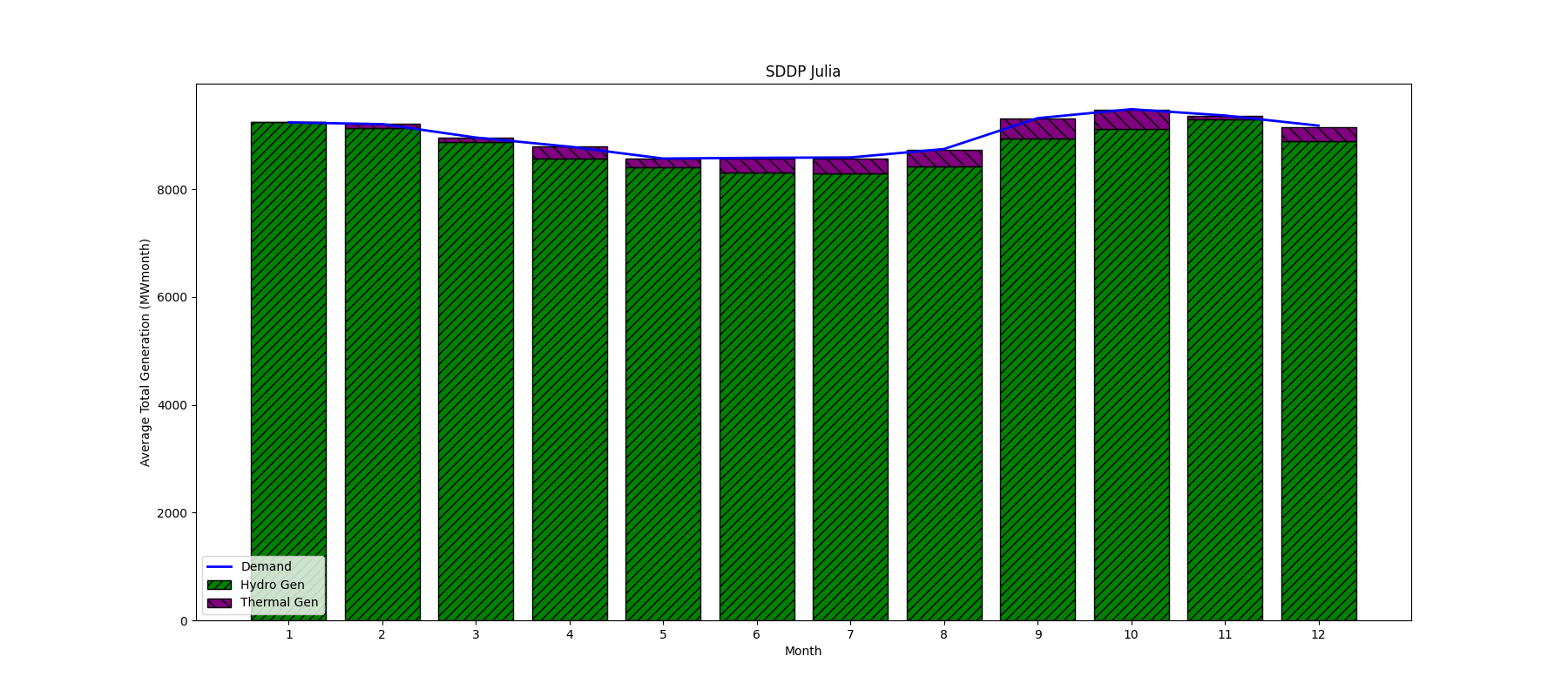}
\caption{Average total generation over 100 simulations for case 1} \label{fig-avgtotalgen-case01}
\end{figure}

We compared the results from our Julia model with those obtained from the {\sc newave} model, as presented in Figures~\ref{fig-zinf-case01}-\ref{fig-avgtotalgen-case01}. Figure~\ref{fig-zinf-case01}, shows the progression of the SDDP lower bound (Zinf) by iteration for both cases. This figure indicates that both SDDP Julia and {\sc newave} converged to the same result, validating the consistency between our model and the {\sc newave} model.

Figure~\ref{fig-cmo-ear-case01} presents the mean value and box plots with the Marginal Operation Cost (MOC, in R\$/MWmonth) and the System Stored Energy (SSE, in MWmonth) for each month and both cases. Again, the results from Julia and {\sc newave} are statistically equivalent, further strengthening the validation of our implementation. Finally, Figure~\ref{fig-avgtotalgen-case01} displays the average hydrothermal generation obtained in Julia over all the simulations. This figure illustrates that the average generation from our model across all 100 simulations meets the demand, which demonstrates the robustness of our approach.

\begin{figure} [hbt]
\centering
\includegraphics[width=1.00\textwidth]{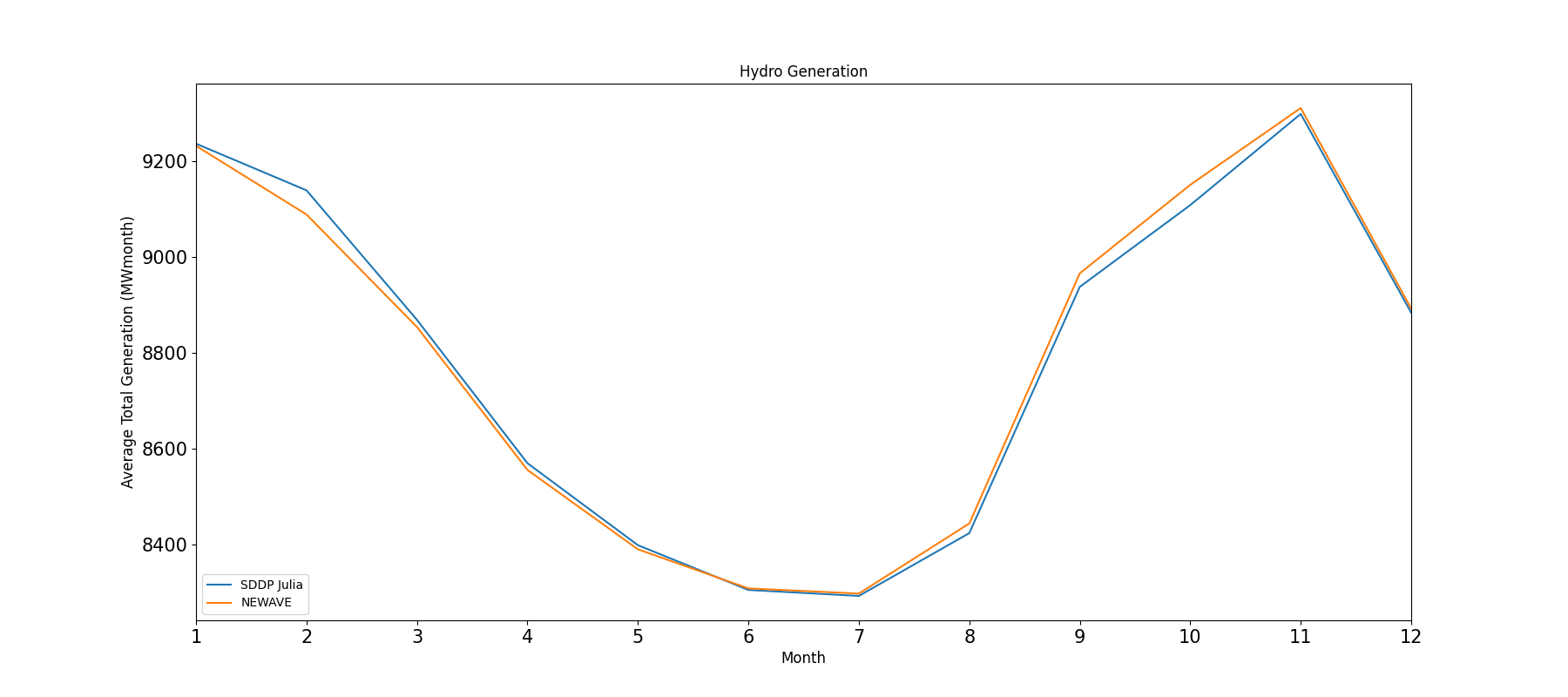}
\includegraphics[width=1.00\textwidth]{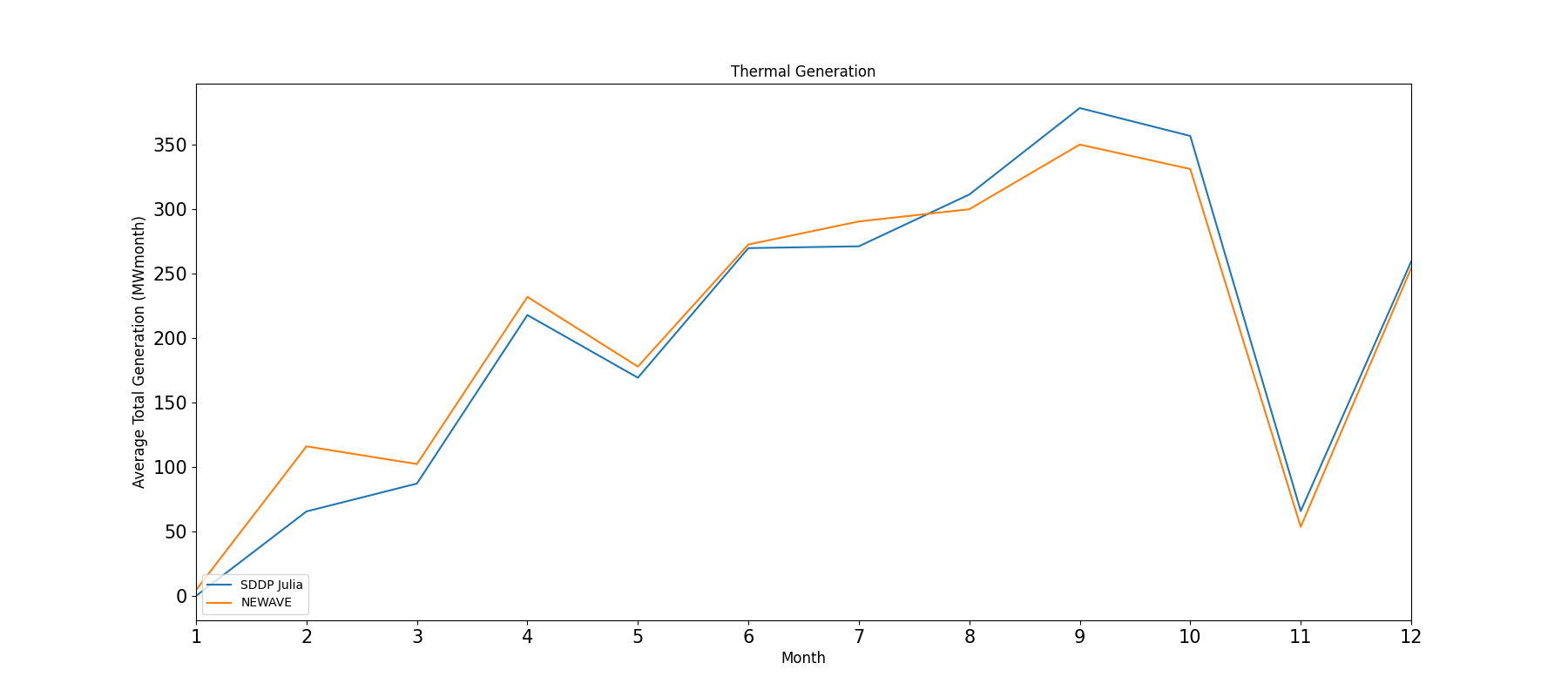}
\caption{Average hydro and thermal total generation over 100 simulations for {\sc sddp.jl} and {\sc newave}.} 
\label{fig-avg-hydrothermal-case01}
\end{figure}

Finally, Figure~\ref{fig-avg-hydrothermal-case01} shows the total average generation for both sources (hydro and thermal), and type of simulations {\sc newave} and {\sc SDDP.jl}. For these graphics, we can ensure that both simulations are consistent.

It should be noted, that the academic version of {\sc newave} used for our benchmark has much less functionalities than the official computational tool. Also, the number of iterations (\# of ite) required by SDDP Julia to reach a bound similar to that of {\sc newave} is higher due to the fact that ten simulations were considered in each forward pass performed by {\sc newave} while a single one was considered in our Julia implementation.

\subsection{Impact of EC}\label{sec:impact_ec}

Following the initial model validation, we simulated the problem incorporating the EC presented in this report, employing the two proposed formulations. To construct the operation policy, we let SDDP run for 400 iterations.

The simulations encompassed different levels of demand and varied values for the penalties of the slack variables added to the outflow limit constraints. These slack variables are a necessary inclusion since the {\sc sddp.jl} package requires the generated subproblems to offer complete recourse. To circumvent the ``end-of-world effect'', as noted in Remark~\ref{rem-abs1}, the optimization horizon extended to $2T= 24$ months, but the analysis horizon remained at $T=12$ months. The simulations performed are listed as follows.

\begin{itemize}
    \item Case 1: three different levels of demand, no EC.
    \item Case 2: three different levels of demand, penalty cost of slack variables equal to 100,  EC  modeled by formulation A.
    \item Case 3: three different levels of demand, penalty cost of slack variables equal to 5000,  EC  modeled by formulation A.
    \item Case 4: three different levels of demand, penalty cost of slack variables equal to 10000,  EC  modeled by formulation A.
    \item Case 5: one level of demand, three penalty cost of slack variables,  EC  modeled by formulation B.
    \item Case 6: one level of demand, no penalty cost for slack variables, no EC; state variables discretized (SDDiP approach) in three different number of samples (10, 50 and 100 samples).
\end{itemize}

\begin{figure} [hbt]
\centering
\includegraphics[width=1.00\textwidth]{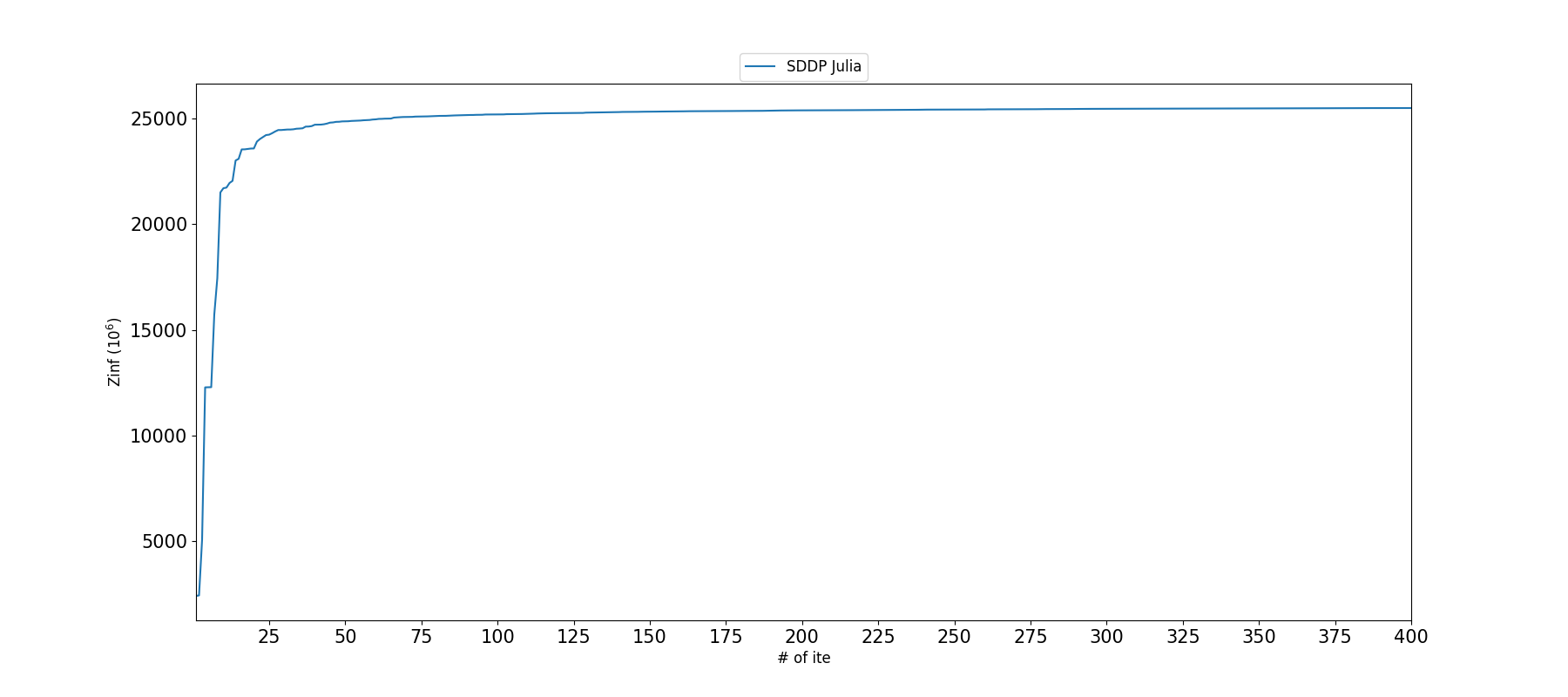}
\caption{Evolution of the lower bound for case 2 with reference demand equal to 9000MW} \label{fig-zinf-case07}
\end{figure}

\begin{figure} [hbt]
\centering
\includegraphics[width=1.00\textwidth]{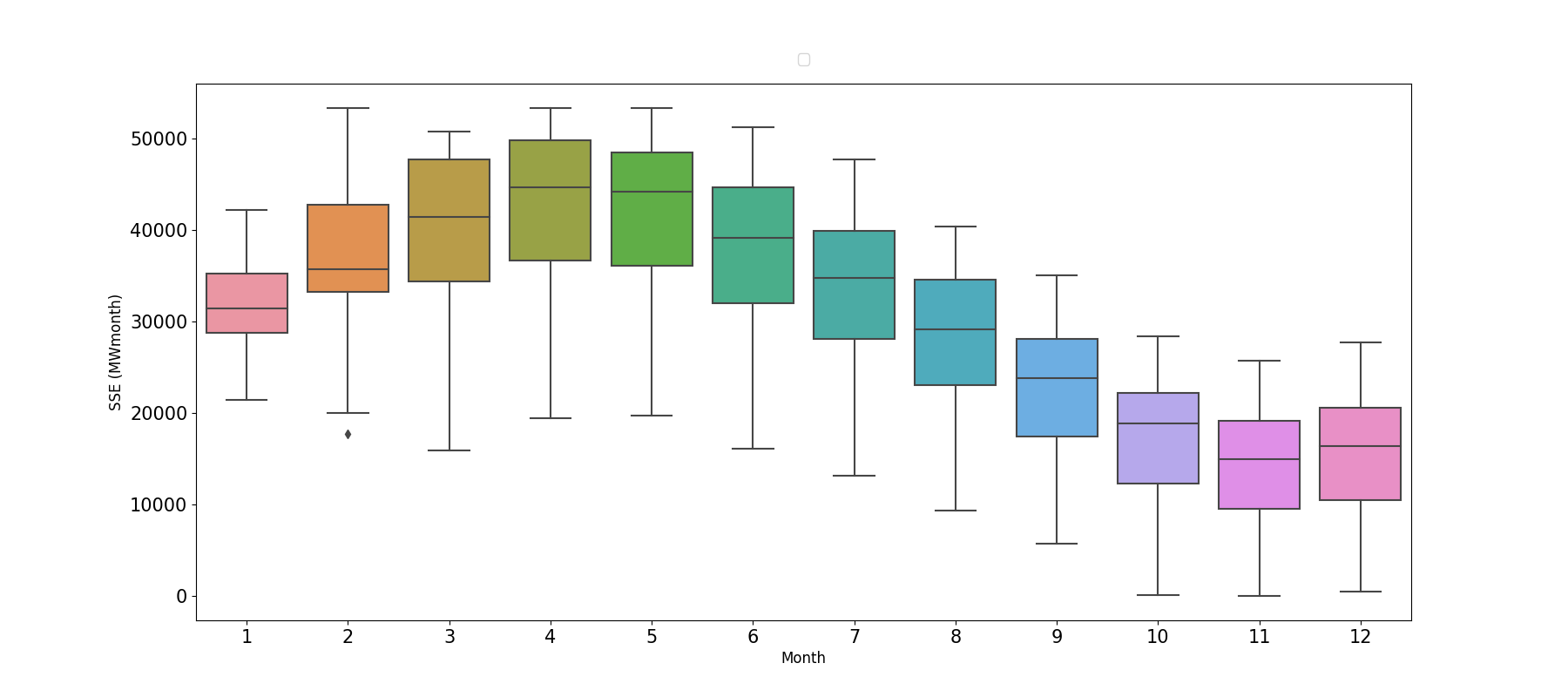}
\caption{System stored energy, mean and boxplots over 100 simulations for case 2.} 
\label{fig-cmo-ear-case07}
\end{figure}

\begin{figure} [hbt]
\centering
\includegraphics[width=1.00\textwidth]{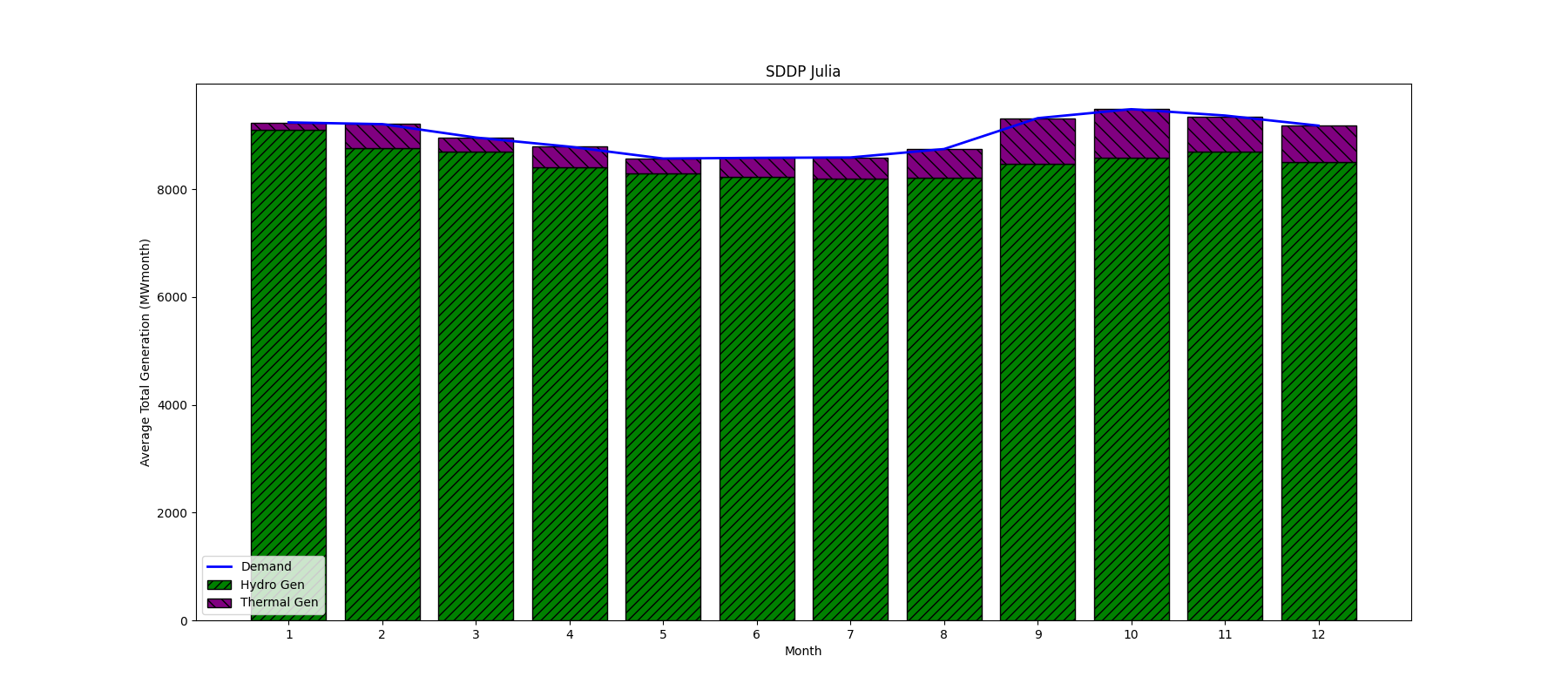}
\caption{Average total generation over 100 simulations for case 2} \label{fig-avgtotalgen-case07}
\end{figure}

The results for Case 2 are depicted in Figures~\ref{fig-zinf-case07}-\ref{fig-avgtotalgen-case07}.

In Figure~\ref{fig-zinf-case07}, we see that {\sc sddp.jl} converges with just over 100 iterations. However, when compared to the results of Section \ref{sec:newave} (without EC), we observe that the lower bound (Zinf) is fifty times larger. This increase results from the penalty applied in the objective function to the slack variables since the {\sc sddp.jl} package requires complete recourse in each subproblem. For this experiment, the penalty value is relatively low, set at 100.

Figure~\ref{fig-cmo-ear-case07}, shows an increase in the stored energy in the system's latter months when compared to case without EC. This rise can be attributed to the imposed EC, which tend to increase water storage in reservoirs.

Figure~\ref{fig-avgtotalgen-case07}, presents the hydrothermal generation over 12 months. The chart illustrates an increase in thermal generation in the final months and, consequently, a decrease in hydroelectric generation. This shift supports the observed trend of higher reservoir storage and underlines the influence of EC on the system's operation.

These results reveal the substantial impact of Environmental Constraints on the system's operation and resource allocation, particularly evident in the increased System Stored Energy (SSE) and shift in hydrothermal generation balance. These observations reinforce the importance of accounting for these factors in long-term hydrothermal planning models.





\begin{figure} [hbt]
\centering
\includegraphics[width=1.00\textwidth]{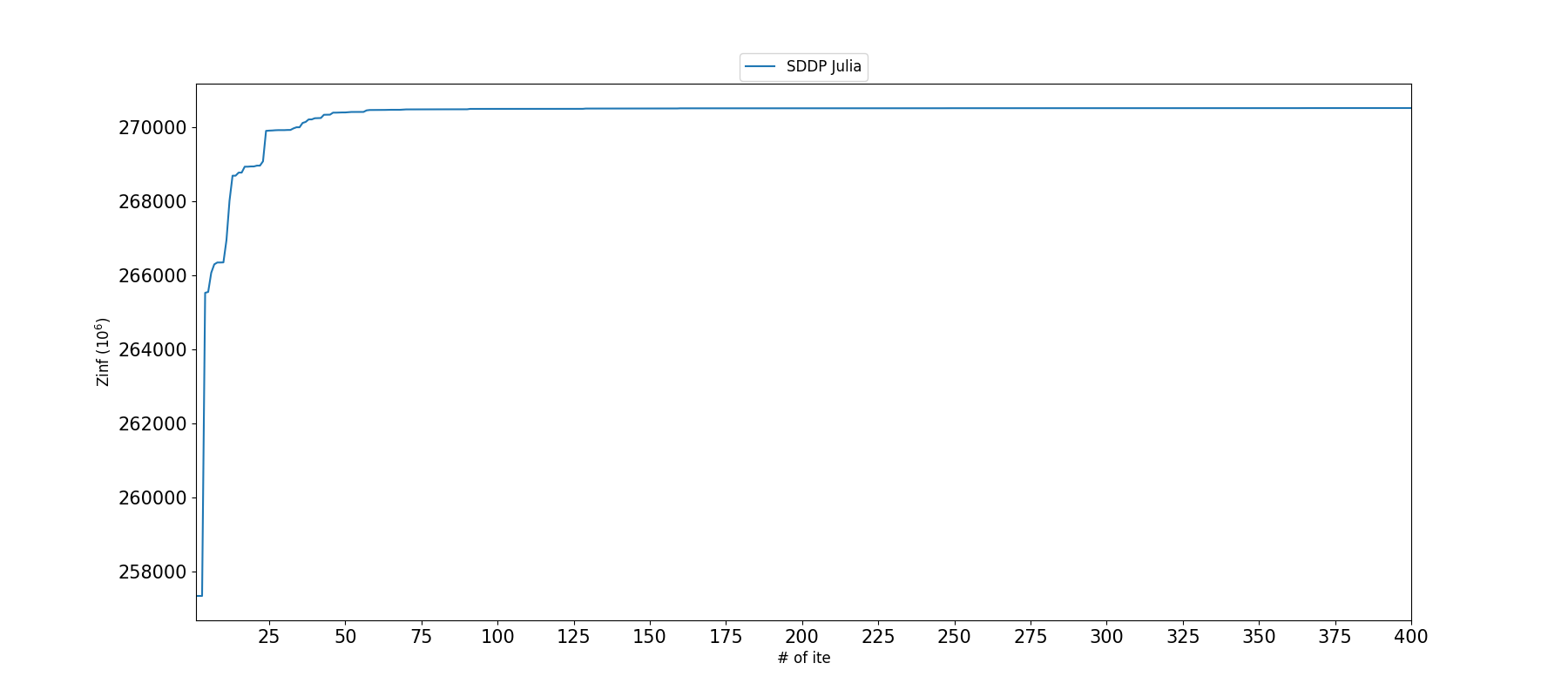}
\caption{Evolution of the lower bound for case 5 with reference demand equal to 9000MW and penalty cost of slack variables equal to 100} \label{fig-zinf-case17}
\end{figure}

\begin{figure} [hbt]
\centering
\includegraphics[width=1.00\textwidth]{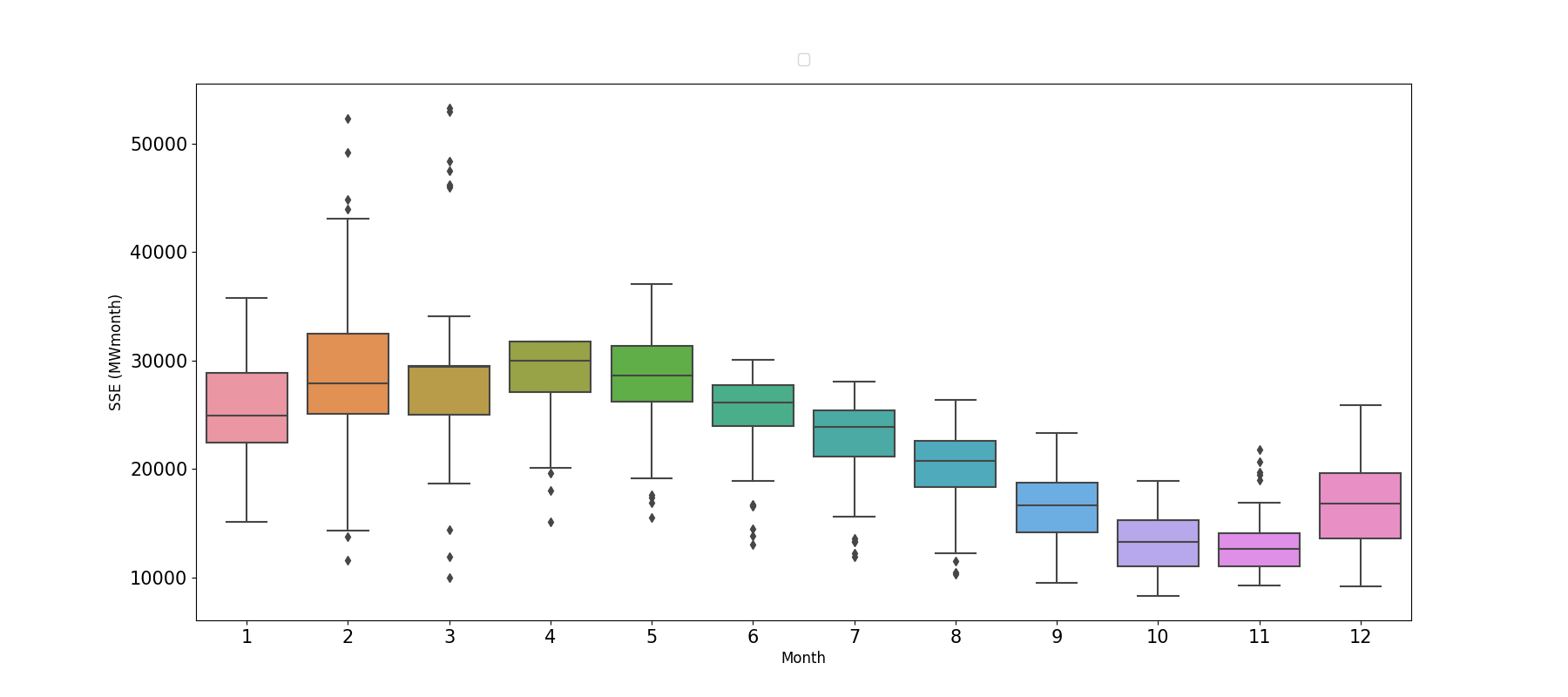}
\caption{System stored energy, mean and boxplots over 100 simulations for case 5 with reference demand equal to 9000MW and penalty cost of slack variables equal to 100} 
\label{fig-cmo-ear-case17}
\end{figure}

\begin{figure} [hbt]
\centering
\includegraphics[width=1.00\textwidth]{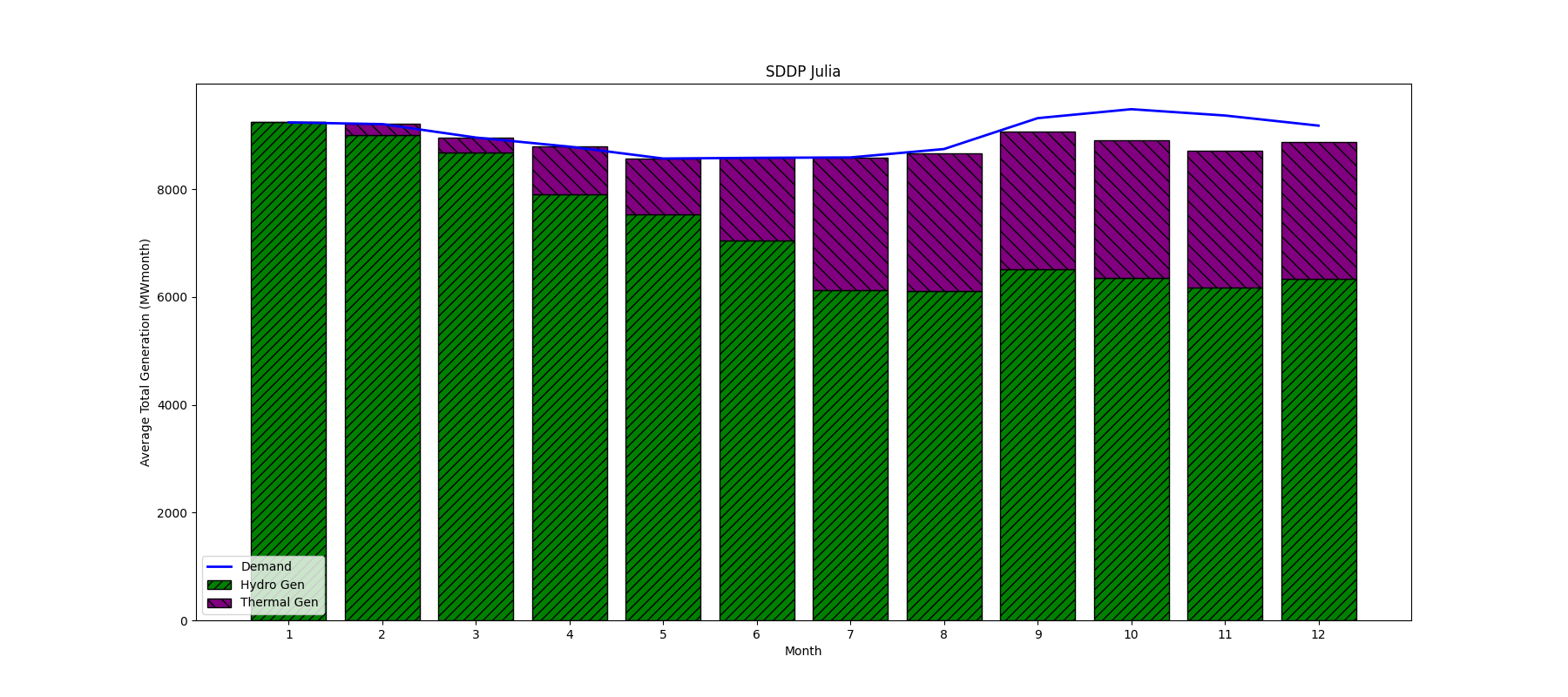}
\caption{Average total generation over 100 simulations for case 5 with reference demand equal to 9000MW and penalty cost of slack variables equal to 100} \label{fig-avgtotalgen-case17}
\end{figure}

We continue our exploration of the impact of environmental constraints by discussing Case 5, which utilizes formulation B to model EC, with a penalty cost of 100 for the slack variables. The results from this case are illustrated in Figures \ref{fig-zinf-case17}-\ref{fig-avgtotalgen-case17}.

In Figure~\ref{fig-zinf-case07}, we observe that the {\sc sddp.jl} algorithm still converges after approximately 100 iterations. However, the lower bound (Zinf) is approximately 11 times higher than Case 2. This disparity can be ascribed to the different EC formulations, which further underscores the critical impact of modeling choices on system operation and cost.

Figure~\ref{fig-cmo-ear-case17} indicates a similar trend to Case 2, with increased water storage, especially in the latter months. Nevertheless, there is a smaller standard deviation compared to Case 2 across the 100 simulations analyzed. This suggests that, for Case 5, there is an even stronger tendency for water to be stored in the reservoirs.

Figure~\ref{fig-avgtotalgen-case07} presents an increase in thermal generation and, consequently, a decrease in hydroelectric generation compared to Case 2. More notably, in the final months, we observe that hydrothermal generation fails to meet demand, resulting in a generation deficit.

The analysis of Case 5 underscores the significant influence of the EC modeling approach on system operation. Compared to the results from Section \ref{sec:newave} (without EC) and Case 2 (with EC modeled by formulation A), we observe larger shifts in system operation, most notably the generation deficit in the final months. This stresses the importance of carefully selecting and calibrating the EC modeling approach, given its substantial influence on system costs and feasibility.



\subsection{SDDiP approach}
\label{sec:sddip}

Finally, we implemented the Stochastic Dual Dynamic Integer Programming (SDDiP) approach, discretizing the state variable of volume into a predefined number of sections. The case explored here corresponds to a model without Environmental Constraints (EC), similar to the one discussed in Section \ref{sec:newave}. We sought to evaluate the potential of this recent approach and its viability for our hydrothermal planning model.

\begin{figure} [hbt]
\centering
\includegraphics[width=1.00\textwidth]{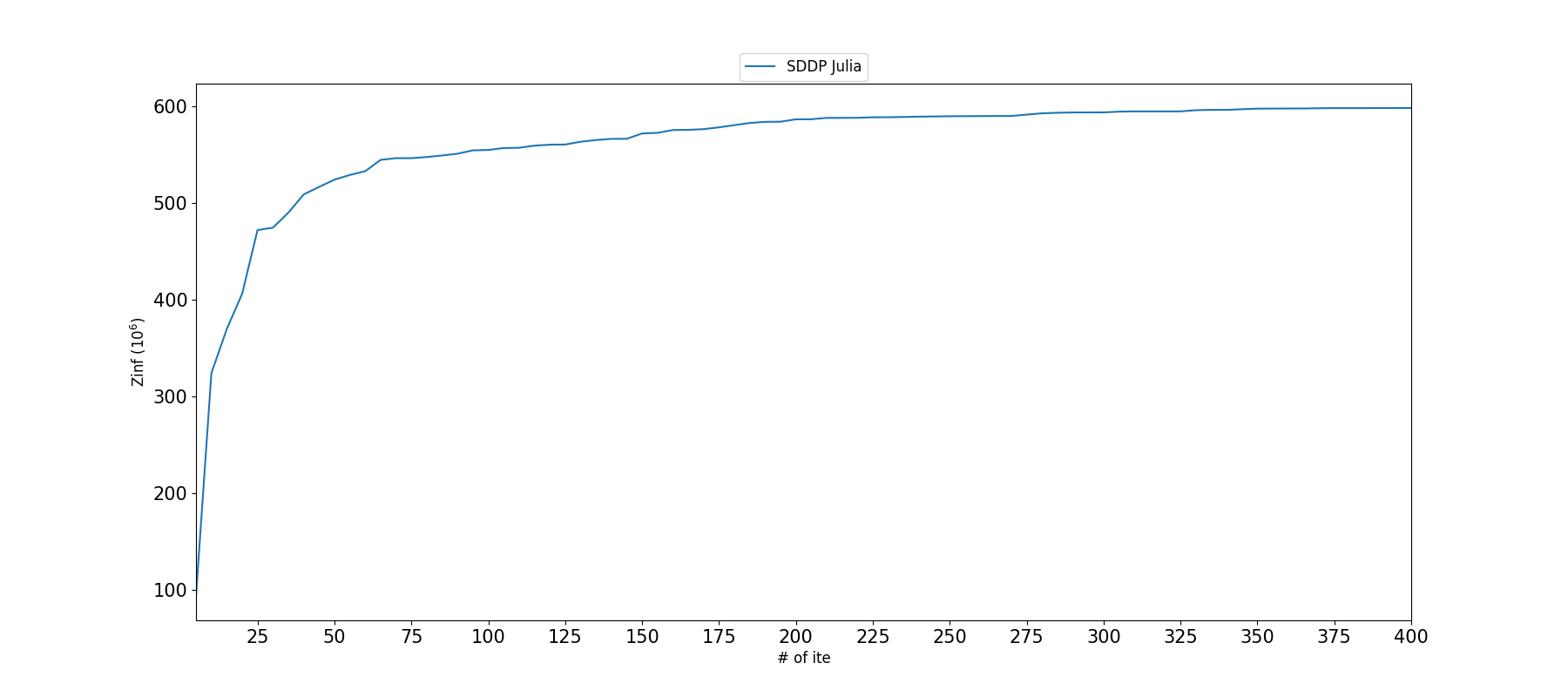}
\caption{Evolution of the lower bound for case 6 with discretization equal to 10 steps} \label{fig-zinf-case21}
\end{figure}

In this experiment, we discretized the volume into ten parts, a relatively low level of discretization. The results obtained from this case are as follows:

Figure~\ref{fig-zinf-case21} shows that the lower bound (Zinf) stabilizes after about 350 iterations. Interestingly, Zinf is slightly higher than the result obtained in Section \ref{sec:newave}. This can be attributed to the approximation error introduced by the low-precision discretization, where the volume is divided into only ten parts.

\begin{figure} [hbt]
\centering
\includegraphics[width=1.00\textwidth]{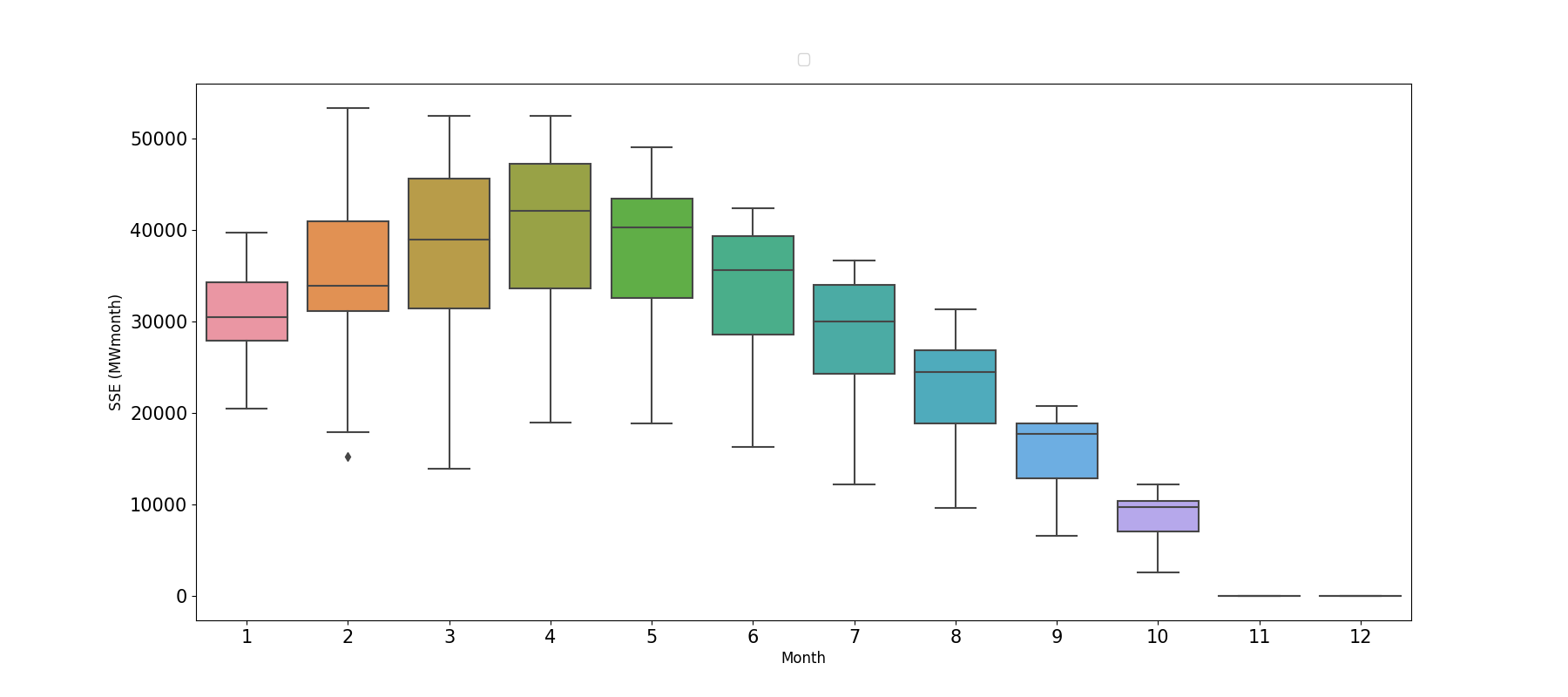}
\caption{System stored energy, mean and boxplots over 100 simulations for case 6 with discretization equal to 10 steps} 
\label{fig-cmo-ear-case21}
\end{figure}

In Figure~\ref{fig-cmo-ear-case21}, we observe higher SSE values compared to Case \ref{sec:newave}. This behavior is also a result of the low-precision discretization strategy implemented in this case.

\begin{figure} [hbt]
\centering
\includegraphics[width=1.00\textwidth]{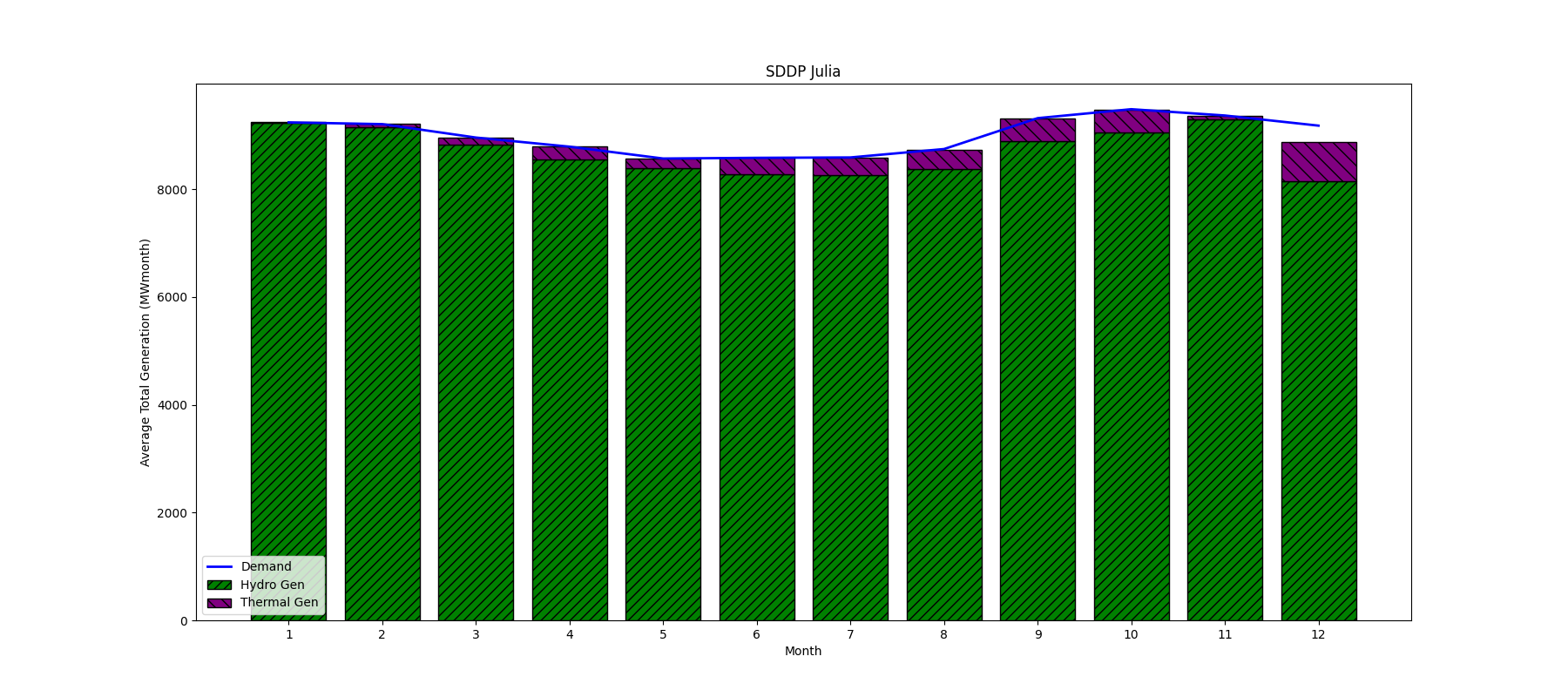}
\caption{Average total generation over 100 simulations for case 6 with discretization equal to 10 steps} 
\label{fig-avgtotalgen-case21}
\end{figure}

In Figure~\ref{fig-avgtotalgen-case21}, we notice a similar behavior to the results of Section 7.1 in terms of hydrothermal generation. However, there is a notable generation deficit in the final month.

The results obtained from the SDDiP approach give us important insights into the trade-offs between computational efficiency and solution accuracy in hydrothermal planning models. While the SDDiP approach converges faster, the low-precision discretization leads to higher lower bounds and system stored energy. Moreover, the approximation error introduced by the low-precision discretization could potentially compromise the model's ability to meet demand, as evidenced by the generation deficit observed in the final month.

Following the investigation into the Stochastic Dual Dynamic Integer Programming (SDDiP) approach, a further simulation was carried out to investigate the impact of higher levels of discretization on the model's results. In this case, the state variable volume was discretized into 100 parts, a significant increase from the previous simulation where we used only 10 parts.

\begin{figure} [hbt]
\centering
\includegraphics[width=1.00\textwidth]{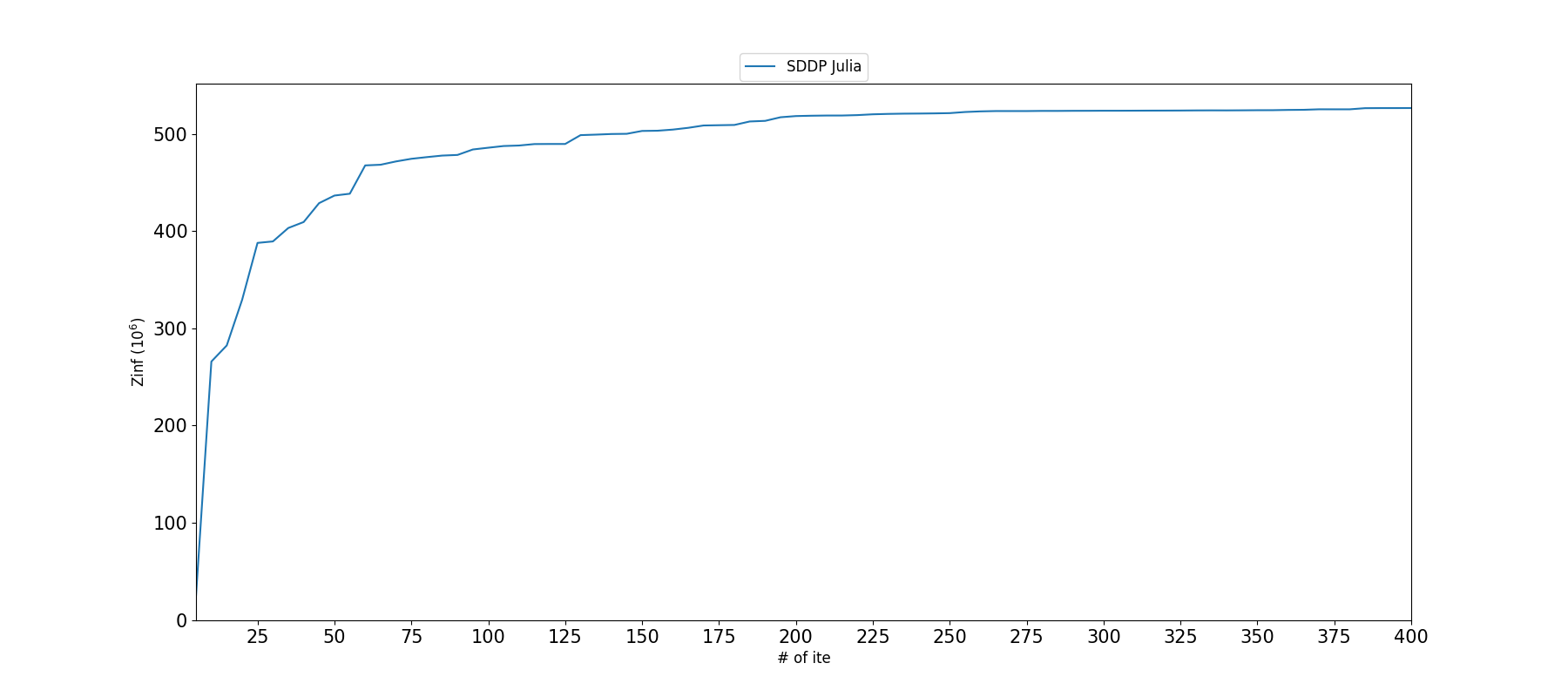}
\caption{Evolution of the lower bound for case 6 with discretization equal to 100 steps} \label{fig-zinf-case23}
\end{figure}

The results of this simulation presented a certain level of consistency with the findings from Section \ref{sec-newave}, particularly in terms of the lower bound (Zinf), as can be seen in Figure \ref{fig-zinf-case23} . This similarity suggests that the increased discretization brought the SDDiP model's lower bound closer to the result from the SDDP model, reinforcing the notion that more precise discretization can improve the model's performance.

\begin{figure} [hbt]
\centering
\includegraphics[width=1.00\textwidth]{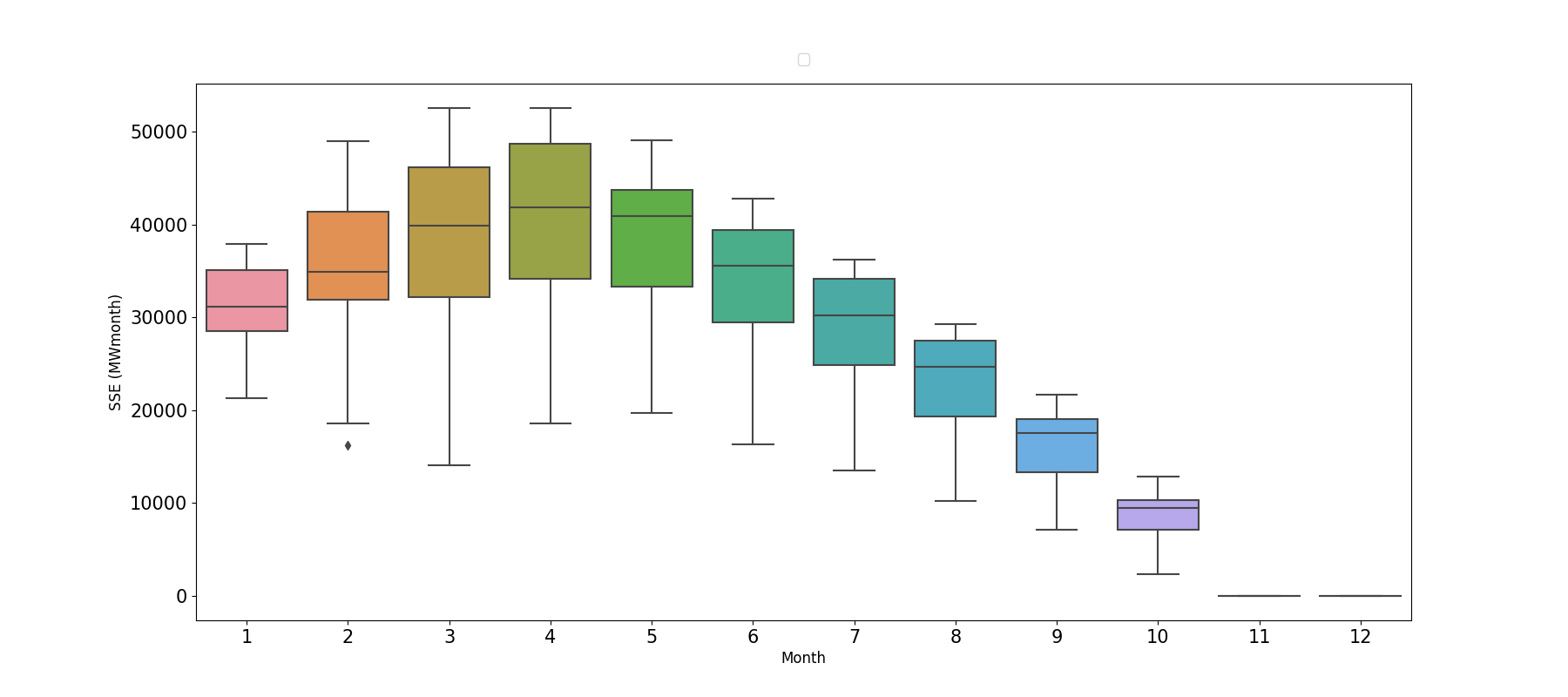}
\caption{System stored energy, mean and boxplots over 100 simulations for case 6 with discretization equal to 100 steps} 
\label{fig-cmo-ear-case23}
\end{figure}

\begin{figure} [hbt]
\centering
\includegraphics[width=1.00\textwidth]{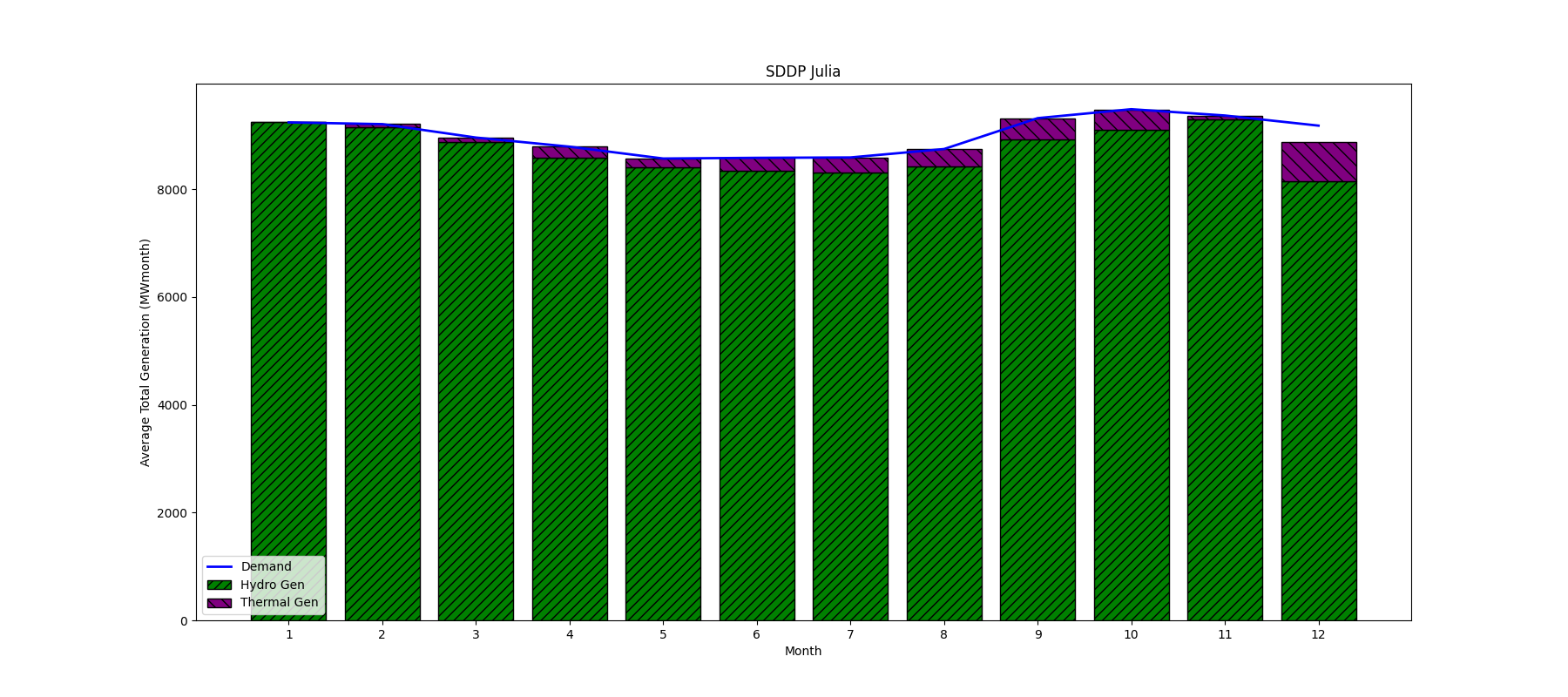}
\caption{Average total generation over 100 simulations for case 6 with discretization equal to 100 steps} 
\label{fig-avgtotalgen-case23}
\end{figure}

However, this was not the case for the other key metrics evaluated. Both the hydrothermal generation (Figure \ref{fig-avgtotalgen-case23}) and system stored energy (Figure \ref{fig-cmo-ear-case23}) metrics demonstrated similar behavior to those obtained with the SDDiP approach at a lower discretization level (10 intervals). This finding suggests that these metrics may not linearly improve with increased discretization, or that a higher level of discretization may be required to attain results comparable to the SDDP approach.

In addition to the effect on the metrics, it is crucial to highlight the computational cost associated with the increased discretization. A larger number of discretization intervals inherently requires higher computational resources, leading to an increase in the model's execution time. This factor should be taken into account when deciding on the level of discretization to use, especially in practical applications where computational resources and execution times are often significant constraints.

These results emphasize the importance of carefully choosing the level of discretization in the SDDiP approach. While a higher level of discretization seems to improve the performance in terms of the lower bound, it may not lead to improvements across all metrics and comes at the cost of increased computational demand.

\section*{Conclusions}\label{sec-conc}

This report delved into the long-term generation scheduling (LTGS) problem within the context of the Brazilian energy system, specifically focusing on the impact of environmental constraints that change the use of hydroelectric reservoirs. By introducing two distinct formulations, A and B, to model these environmental constraints, we demonstrated how the complexity of constraints can affect the final solution of the optimization problem.

Each formulation brought along its own set of advantages and disadvantages regarding the representation of mathematical functions, usage of binary variables, and the quantity of constraints. These factors directly impacted the number of slack variables employed due to the complete recourse requirement of the {\sc SDDP.jl} toolbox, and also influenced the overall cost of the problem due to penalties associated with the slack variables.

Through the experiments conducted, we observed that the inclusion of environmental constraints affects the lower bounds, system stored energy, and the hydrothermal generation. Our study showed that the choice of formulation and penalty parameter for slack variables significantly influences the LTGS problem's solution.

We found that an increase in the penalty parameter strengthens the environmental constraints, leading to a reduction in outflows, an increase in reservoir volumes, and subsequently, an increase in system stored energy. As the penalty parameter escalates, violating the outflow constraints becomes economically infeasible, leading to a decrease in hydroelectric generation. To counterbalance the drop in hydroelectric generation, thermal generation increases to meet the energy demand.

Furthermore, our investigation into the Stochastic Dual Dynamic Integer Programming (SDDiP) approach highlighted the need for careful consideration of the level of discretization. While the SDDiP strategy shows promise for hydrothermal planning models, its effectiveness is influenced by the precision of discretization.

Several avenues for future research emerged from our study:
\begin{itemize}
    \item Exploring the impact of the discrepancy in the number of forward pass simulations between {\sc NEWAVE} and our Julia implementation could provide valuable insights into computational efficiency, solution quality, and overall model performance.

    \item Further exploration into the impact of the penalties used for the slack variables, and the development of additional modeling strategies to handle the complexities introduced by environmental constraints could improve model accuracy and robustness.

    \item The differences between Formulation A and Formulation B under environmental constraints warrant further investigation to fully understand their implications for hydrothermal planning models.
In particular, it is interesting to check more carefully the impact of the minimum inflow relaxation when the reservoir level is very low.

    \item Studying the impact of higher precision discretization on the accuracy, robustness, and computational efficiency of the SDDiP approach could yield significant benefits.

    \item The performance of the SDDiP approach under different case conditions, such as models incorporating Environmental Constraints, could give a broader understanding of its potential benefits and limitations.
\end{itemize}

In conclusion, this report underscored the criticality of adequately modeling environmental constraints in long-term hydrothermal planning. As our understanding of these constraints deepens, more nuanced modeling strategies will need to be developed to better accommodate these complexities, ensuring that our models remain robust, accurate, and capable of guiding energy system planning effectively.

\printbibliography
\begin{appendix}
\section{Data of test system}\label{sec:test-system-data}

The hydro configuration of the test system in Figure \ref{fig-hydro-conf} has
data given in Table \ref{tab:hydro-data-testsystem}.

\begin{table}
\centering
\begin{tabular}{cccccc}
\hline
Plant        & \begin{tabular}[c]{@{}c@{}}Maximum\\ generation\\  (MWmonth)\end{tabular} & \begin{tabular}[c]{@{}c@{}}Useful\\ volume (hm$^3$)\end{tabular} & \begin{tabular}[c]{@{}c@{}}Productivity\\ (MWmonth/hm$^3$)\end{tabular} & \begin{tabular}[c]{@{}c@{}}Maximum\\ turbined outflow\\ (hm$^3$)\end{tabular} & \begin{tabular}[c]{@{}c@{}}Initial\\ volume\\ (hm$^3$)\end{tabular} \\ \hline
RETIRO BAIXO & 82.74256                                                                  & 40.87                                                         & 0.12723                                                         & 662.20071                                                                  & 27.38                                                            \\ 
TRES MARIAS  & 362.88700                                                                 & 15278.00                                                      & 0.16576                                                         & 2112.89078                                                                 & 7654.28                                                          \\
QUEIMADO     & 100.65465                                                                 & 461.75                                                        & 0.62372                                                         & 166.44477                                                                  & 258.12                                                           \\
SOBRADINHO   & 1008.87562                                                                & 28669.00                                                      & 0.09312                                                         & 10690.01959                                                                & 14821.87                                                         \\ 
ITAPARICA    & 1422.38577                                                                & 3548.00                                                       & 0.16826                                                         & 8263.97182                                                                 & 2117.09                                                          \\
COMP PAF-MOX & 3883.02034                                                                & 0.00                                                          & 0.38833                                                         & 9999.35762                                                                 & 0.00                                                             \\
XINGO        & 3079.08318                                                                & 0.00                                                          & 0.40992                                                         & 7511.35249                                                                 & 0.00                                                             \\ \hline
\end{tabular}
\caption{Data from hydros of the test system.}\label{tab:hydro-data-testsystem}
\end{table}

The thermal data of the test system is presented in Table \ref{tab:therm-data-testsystem}.

\begin{table}
\centering
\begin{tabular}{ccc}
\hline
\multicolumn{1}{c}{Plant} &
\multicolumn{1}{c}{\begin{tabular}[c]{@{}c@{}}Maximum \\ generation\\
(MW)\end{tabular}} & \multicolumn{1}{c}{\begin{tabular}[c]{@{}c@{}}$c_g$\\
(constant in time)\\ (R\$/MWh)\end{tabular}} \\ \hline
GOIANIA II                  & 136                                                                                       & 1927                                                                          \\
JUIZ DE FORA                & 87                                                                                        & 523                                                                           \\
PALMEIRAS GO                & 140                                                                                       & 1493                                                                          \\
TERMOMACAE                  & 929                                                                                       & 880                                                                           \\
VIANA                       & 437                                                                                       & 1257                                                                          \\
XAVANTES                    & 54                                                                                        & 2632                                                                          \\
CAMPINA GDE                 & 169                                                                                       & 1257                                                                          \\
PETROLINA                   & 594                                                                                       & 2014  \\ \hline                                     
\end{tabular}
\caption{Data from thermal plants of the test system.}\label{tab:therm-data-testsystem}
\end{table}

All hydrothermal plants are connected, resulting in a single bus system with:
\begin{itemize}
    \item 1 deficit level, with associated cost equal to 7643.82 \$/MWh
    \item 1 load level
    \item The load, in MWmonth, is given in Table \ref{tab:load-data-testsystem}
\end{itemize}

\begin{table}[]
\centering
\begin{tabular}{c|ccc}
\hline
\multirow{2}{*}{Month} & \multicolumn{3}{c}{Load level (mean)}                           \\ \cline{2-4} 
                       & \multicolumn{1}{c|}{8,000} & \multicolumn{1}{c|}{9,000} & 10,000 \\ \hline
1                      & \multicolumn{1}{c|}{8,211} & \multicolumn{1}{c|}{9,237} & 10,263 \\ \hline
2                      & \multicolumn{1}{c|}{8,182} & \multicolumn{1}{c|}{9,205} & 10,227 \\ \hline
3                      & \multicolumn{1}{c|}{7,961} & \multicolumn{1}{c|}{8,956} & 9,951  \\ \hline
4                      & \multicolumn{1}{c|}{7,811} & \multicolumn{1}{c|}{8,787} & 9,763  \\ \hline
5                      & \multicolumn{1}{c|}{7,615} & \multicolumn{1}{c|}{8,567} & 9,519  \\ \hline
6                      & \multicolumn{1}{c|}{7,626} & \multicolumn{1}{c|}{8,580} & 9,533  \\ \hline
7                      & \multicolumn{1}{c|}{7,633} & \multicolumn{1}{c|}{8,587} & 9,541  \\ \hline
8                      & \multicolumn{1}{c|}{7,772} & \multicolumn{1}{c|}{8,743} & 9,715  \\ \hline
9                      & \multicolumn{1}{c|}{8,281} & \multicolumn{1}{c|}{9,316} & 10,351 \\ \hline
10                     & \multicolumn{1}{c|}{8,429} & \multicolumn{1}{c|}{9,482} & 10,536 \\ \hline
11                     & \multicolumn{1}{c|}{8,324} & \multicolumn{1}{c|}{9,365} & 10,405 \\ \hline
12                     & \multicolumn{1}{c|}{8,157} & \multicolumn{1}{c|}{9,177} & 10,196 \\ \hline
\end{tabular}
\caption{Data of demand of the test system.}\label{tab:load-data-testsystem}
\end{table}

\section{Simulation Results}

\begin{table}
\centering
\begin{tabular}{|cccccc|}
\hline
Simulation & Load level & Penalty & Planning horizon & Zinf & Runtime\\
case & (mean) & cost (\$) & (months) & (\$) & (s) \\ \hline
1 &  9,000 &      0 & 24 & $6.01 \times 10^9$    & 32.89  \\ \hline
1 &  8,000 &      0 & 24 & $1.25 \times 10^8$    & 37.38  \\ \hline
1 & 10,000 &      0 & 24 & $2.39 \times 10^{10}$ & 32.38  \\ \hline
2 &  9,000 &    100 & 24 & $2.55 \times 10^{10}$ & 90.61  \\ \hline
2 &  8,000 &    100 & 24 & $1.69 \times 10^{10}$ & 98.77  \\ \hline
2 & 10,000 &    100 & 24 & $4.41 \times 10^{10}$ & 97.90  \\ \hline
3 &  9,000 &  5,000 & 24 & $7.45 \times 10^{11}$ & 84.54  \\ \hline
3 &  8,000 &  5,000 & 24 & $6.93 \times 10^{11}$ & 100.95 \\ \hline
3 & 10,000 &  5,000 & 24 & $8.11 \times 10^{11}$ & 77.59  \\ \hline
4 &  9,000 & 10,000 & 24 & $1.38 \times 10^{12}$ & 73.89  \\ \hline
4 &  8,000 & 10,000 & 24 & $1.32 \times 10^{12}$ & 71.36  \\ \hline
4 & 10,000 & 10,000 & 24 & $1.44 \times 10^{12}$ & 87.46  \\ \hline
5 &  9,000 &      0 & 24 & $1.02 \times 10^{10}$ & 43.08  \\ \hline
5 &  9,000 &    100 & 24 & $2.71 \times 10^{11}$ & 40.56  \\ \hline
5 &  9,000 & 10,000 & 24 & $5.0  \times 10^{11}$ & 42.95  \\ \hline
6 &  9,000 &      0 & 24 & $6.07 \times 10^8$    & 150.17 \\ \hline
6 &  9,000 &      0 & 24 & $5.33 \times 10^8$    & 813.43 \\ \hline
6 &  9,000 &      0 & 24 & $5.33 \times 10^8$    & 1,141.38 \\ \hline
\end{tabular}
\caption{SDDP running times.}\label{tab:sddp-run-times}
\end{table}

\end{appendix}

\end{document}